\documentclass[11pt]{article}
\usepackage{mathptmx}
\usepackage{graphicx, amssymb, amsthm}
\usepackage[papersize={5.5in, 8.5in}, top=0.35in, bottom=0.5in, left=0.35in, right=0.35in]{geometry}

\usepackage{amsmath}
\usepackage{color}
\usepackage{tikz}
\usepackage{hyperref}

\date{}

\theoremstyle{definition}
\newtheorem{df}{Definition} [section]

\theoremstyle{plain}
\newtheorem{thm}[df]{Theorem}

\newtheorem{problem}[df]{Problem}

\usepackage{algorithm}
\usepackage{algpseudocode}
\usepackage{tikz}
\usepackage{hyperref}

\definecolor{ududff}{rgb}{0.30196078431372547,0.30196078431372547,1.}
\definecolor{ffqqqq}{rgb}{1.,0.,0.}
\definecolor{ffffff}{rgb}{1.,1.,1.}
\definecolor{ccqqww}{rgb}{0.8,0.,0.4}
\definecolor{ttzzqq}{rgb}{0.2,0.6,0.}
\definecolor{ffqqqq}{rgb}{1.,0.,0.}
\definecolor{qqqqff}{rgb}{0.,0.,1.}
\definecolor{ududff}{rgb}{0.30196078431372547,0.30196078431372547,1.}
\definecolor{ffffff}{rgb}{1.,1.,1.}
\definecolor{ududff}{rgb}{0.30196078431372547,0.30196078431372547,1.}
\definecolor{ffqqqq}{rgb}{1.,0.,0.}
\definecolor{ffffff}{rgb}{1.,1.,1.}
\definecolor{ududff}{rgb}{0.30196078431372547,0.30196078431372547,1.}
\definecolor{ffqqqq}{rgb}{1.,0.,0.}
\definecolor{ffffff}{rgb}{1.,1.,1.}

\definecolor{yqyqyq}{rgb}{0.5019607843137255,0.5019607843137255,0.5019607843137255}
\definecolor{uuuuuu}{rgb}{0.26666666666666666,0.26666666666666666,0.26666666666666666}
\definecolor{zzttqq}{rgb}{0.6,0.2,0.}
\definecolor{zzwwqq}{rgb}{0.6,0.4,0.}
\definecolor{ttffcc}{rgb}{0.2,1.,0.8}
\definecolor{ffqqff}{rgb}{1.,0.,1.}
\definecolor{qqwuqq}{rgb}{0.,0.39215686274509803,0.}
\definecolor{qqzzqq}{rgb}{0.,0.6,0.}
\definecolor{ffqqqq}{rgb}{1.,0.,0.}
\definecolor{ffffqq}{rgb}{1.,1.,0.}
\definecolor{qqffff}{rgb}{0.,1.,1.}
\definecolor{ttzzqq}{rgb}{0.2,0.6,0.}
\definecolor{ffxfqq}{rgb}{1.,0.4980392156862745,0.}
\definecolor{qqqqff}{rgb}{0.,0.,1.}
\definecolor{ffffff}{rgb}{1.,1.,1.}
\definecolor{ffffff}{rgb}{1.,1.,1.}
\definecolor{ududff}{rgb}{0.30196078431372547,0.30196078431372547,1.}

\title{The chromatic number of the Minkowski plane - the regular polygon case.}
\author{%
Geoffrey Exoo\\
Department of Mathematics and Computer Science\\
Indiana State University\\
Terre Haute, IN 47809,
\texttt{ge@cs.indstate.edu}
\and
David Fisher\\
Hofstra University\\
Hempstead, NY 11549,
\texttt{dfisher1@pride.hofstra.edu}
\and
Dan Ismailescu\\
Mathematics Department\\
Hofstra University\\
Hempstead, NY 11549,
\texttt{dan.p.ismailescu@hofstra.edu}
}

\begin{document}

\maketitle
\thispagestyle{empty}
\pagestyle{empty}

\begin{abstract}
The Hadwiger-Nelson problem asks for the minimum number of colors, so that each point of the plane can be assigned a single color with the property that no two points unit-distance apart are identically colored. It is now  known that the answer is $5$, $6$, or $7$, Here we consider the problem in the context of Minkowski planes, where the unit circle is a regular polygon with $8$, $10$, or $12$ vertices. We prove that in each of these cases, one also needs at least five colors.
\end{abstract}

\section{\bf Introduction}
In 1950, Edward Nelson raised the problem of determining the minimum number of colors that are needed to color the points of the Euclidean plane so that no two points unit distance apart are assigned the same color. This number is referred to as the {\it chromatic number of the plane}, and is denoted by $\chi(\mathbf{R}^2)$. The bounds $4\le \chi(\mathbf{R}^2)\le 7$ were easily proved shortly afterwards (see \cite{hadwiger, MM}).

There was no progress on the problem until 2018 when de Grey \cite{degrey}, and independently, the authors of this note \cite{exooismailescuchi5} were able to construct the first $5$-chromatic unit-distance graphs. A flurry of activity ensued, pursuing two main goals. On one hand, progressively smaller 5-chromatic unit-distance graphs were found by Heule \cite{heule} and subsequently by Parts \cite{parts1}; the current record is a graph of order $510$. On the other hand, since all these constructions rely on computer assistance it may be desirable to develop techniques that would make a traditional verification possible. A first contribution in this direction is due to Parts \cite{parts2}.
A well documented and entertaining history of the problem is presented by Soifer in his monograph \cite{soifer}. Several variations of the original problem have been recently considered. We direct the interested reader to the Polymath16 research threads \cite{polymath16} for the most current developments.

One interesting avenue of research is to try to extend these results to other metrics. Chilakamarri \cite{chilakamarri} considered the Minkowski plane scenario; we briefly mention his findings below.

Let $C$ be a convex closed curve, symmetric with respect to the origin $\mathbf{o}$ in the Euclidean plane $\mathbf{R}^2$. Let $(\mathbf{R}^2,C)$ denote the Minkowski metric space where $C$ serves as a unit circle. Two points $\mathbf{x}$ and $\mathbf{y}$ in $\mathbf{R}^2$ are said to be Minkowski distance $1$ apart if there exists a point $\mathbf{p} \in C$ such that $\mathbf{x}-\mathbf{y}=\mathbf{p}$. We refer the reader to the survey of Martini et al. \cite{martini} for more information on Minkowski spaces.

The \emph{chromatic number} of this plane, denoted $\chi(\mathbf{R}^2,C)$, is the least number of colors needed to color the points of $\mathbf{R}^2$ such that no two points Minkowski unit distance apart are identically colored. If $C=S^1$, the unit circle with the traditional Euclidean metric, then $\chi(\mathbf{R}^2,S^1)=\chi(\mathbf{R}^2)$. As mentioned earlier, it is now known that $5\le \chi(\mathbf{R}^2,S^1)\le 7$.

By extending the arguments of Moser and Moser \cite{MM} and the construction of Hadwiger \cite{hadwiger} to arbitrary Minkowski planes, Chilakamarri \cite{chilakamarri} proved that for any centrally symmetric convex curve $C$, we have $4\le \chi(\mathbf{R}^2, C)\le 7$. Moreover, he proved that if $C$ is a parallelogram or a  hexagon then $\chi(\mathbf{R}^2, C)=4$.

To the best of our knowledge, these are the sole instances of $C$ for which the exact value the chromatic number of the Minkowski plane is known. Moreover, with the exception of $S^1$ there is no known example of $C$ for which one can prove that $\chi(\mathbf{R}^2,C)\ge 5$.

It is compelling to investigate the problem when the Minkowski unit circle is a centrally symmetric polygon with $\ge 8$ vertices. In fact, one of Chilakamarri's questions is whether it is true that $\chi(\mathbf{R}^2, P_8)=4$ for every centrally symmetric octagon $P_8$. We answer his question in the negative.

In the following sections we prove that $\chi(\mathbf{R}^2,C)\ge 5$ when $C$ is a regular octagon, a regular decagon, or a regular dodecagon. It turns out that it is easier to construct $5$-chromatic unit distance graphs in any of these metrics than it is to achieve the same result in the regular Euclidean metric setting; our best constructions are graphs of order $120, 121$, and $295$, respectively.

In each case, we follow the same approach which we briefly describe below.

First, we search for a relatively small set of Minkowski unit vectors, $U$, such that the unit distance graph with vertex set $U+U=\{\mathbf{x}+\mathbf{y}\,|\, \mathbf{x},\mathbf{y}\in U\}$ and edge set $\{\{\mathbf{u},\mathbf{v}\}\,|\,\mathbf{u},\mathbf{v}\in U+U \,\, \text{and}\,\,\mathbf{u}-\mathbf{v}\in U \}$ has high edge density.

We call $U$ to be the set of \emph{generating} vectors of the unit-distance graph. In each of the three cases we consider,  this graph turns out to be $4$-chromatic.

Second, we look at pairs of points $\mathbf{u},\mathbf{v}\in U+U$ with the property that the difference  $\mathbf{u}-\mathbf{v}$
is a valid Minkowski unit vector \emph{even though $\mathbf{u}-\mathbf{v}$ does not belong to the set $U$!} We thus discover an additional set of unit vectors, $V$, which we refer to as the set of \emph{accidental} unit vectors.

 We then investigate whether the additional unit edges produced by vectors in $V$ affect the chromatic number of the unit-distance graph with vertex set $U+U$. In each of the three cases, the resulting graph is $5$-chromatic.

Finally, in each case we attempt to identify a $5$-chromatic subgraph of order as small as possible. Details are presented in the subsequent sections.

\section{\bf $C=$ regular octagon }

Consider the regular octagon with vertices
\begin{align*}
\mathbf{a}_1&=(4,0), \mathbf{a}_2=(2\sqrt{2},2\sqrt{2}), \mathbf{a}_3=(0,4),\mathbf{a}_4=(-2\sqrt{2},2\sqrt{2}),\\
-\mathbf{a}_1&=(-4,0), -\mathbf{a}_2=(-2\sqrt{2},-2\sqrt{2}), -\mathbf{a}_3=(0,-4),-\mathbf{a}_4=(2\sqrt{2},-2\sqrt{2}).
\end{align*}

Throughout this section we will only deal with points whose coordinates are of the form $(a+b\sqrt{2}, c+d\sqrt{2})$ with $a, b, c, d$ integers.
In order to avoid computations involving $\sqrt{2}$ we will encode such a point in the format $(a,b,c,d)$.
With this convention, the vertices of the regular octagon become
\begin{align*}
\mathbf{a}_1&=(4,0,0,0), \mathbf{a}_2=(0,2,0,2), \mathbf{a}_3=(0,0,4,0),\mathbf{a}_4=(0,-2,0,2),\\
-\mathbf{a}_1&=(-4,0,0,0), -\mathbf{a}_2=(0,-2,0,-2), -\mathbf{a}_3=(0,0,-4,0),-\mathbf{a}_4=(0,2,0,-2).
\end{align*}

For a given point $\mathbf{p}_1=(a,b,c,d)$, let $\mathbf{p_2}$, $\mathbf{p_3}$ and $\mathbf{p}_4$ be the images of $\mathbf{p}_1$ after  counterclockwise rotations of $\pi/4, 2\pi/4$ and $3\pi/4$, respectively.

It is easy to check that
\begin{align}\label{octagonorbit1}
&\mathbf{p}_1=(a,b,c,d), \mathbf{p_2}=\left(b-d,\frac{a-c}{2},b+d, \frac{a+c}{2}\right),\notag\\
&\mathbf{p_3}=(-c,-d,a,b), \mathbf{p}_4=\left(-b-d,-\frac{a+c}{2},b-d,\frac{a-c}{2} \right).
\end{align}

For a given point $\mathbf{p}_1$, we denote by $\langle \mathbf{p}_1\rangle$ the set of vertices of the regular octagon centered at the origin one of whose vertices is $\mathbf{p}_1$.
\begin{equation}\label{octagonorbit2}
\langle \mathbf{p}_1\rangle :=\{\mathbf{p}_1, \mathbf{p}_2, \mathbf{p}_3, \mathbf{p}_4, -\mathbf{p}_1, -\mathbf{p}_2, -\mathbf{p}_3, -\mathbf{p}_4\}
\end{equation}
We refer to $\langle \mathbf{p}_1\rangle$ as the \emph{octagonal orbit} of $\mathbf{p}_1$.

Obviously, $\langle \mathbf{p}_1\rangle=\langle \mathbf{p}_2\rangle=\langle \mathbf{p}_3\rangle=\langle \mathbf{p}_4\rangle=\langle -\mathbf{p}_1\rangle=\langle -\mathbf{p}_2\rangle=\langle -\mathbf{p}_3\rangle=\langle -\mathbf{p}_4\rangle$.

Chilakamarri proved that the Moser spindle can be embedded as a unit-distance graph in any Minkowski metric. We will first present such an embedding in the regular octagon metric.

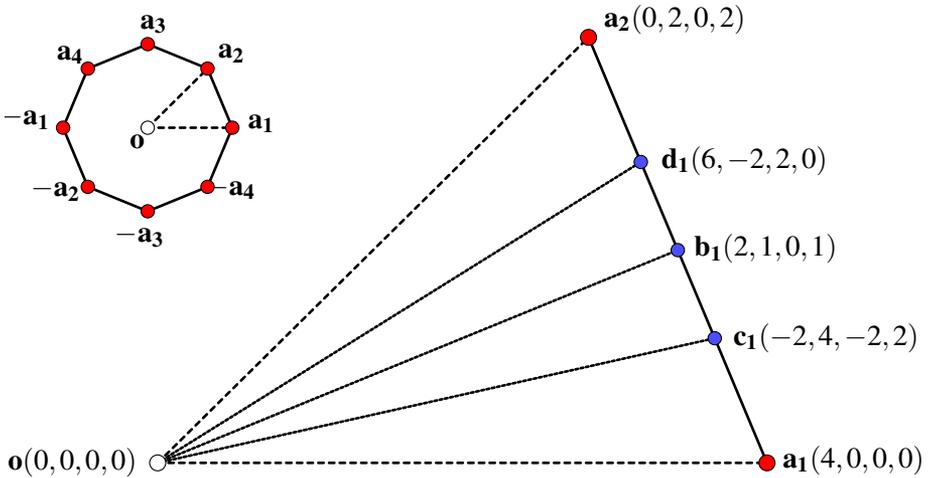
\begin{figure}[ht]
\centering
\begin{tikzpicture}[line cap=round,line join=round,x=1.0cm,y=1.0cm]
\clip(-2.0066666666666677,-0.3177777777777797) rectangle (10.04666666666665,6.588888888888886);
\draw (5.726666666666656,6.224444444444441) node[anchor=north west] {$\mathbf{a_2}(0,2,0,2)$};
\draw (8.07333333333332,0.3666666666666646) node[anchor=north west] {$\mathbf{a_1} (4,0,0,0)$};
\draw (7.424444444444432,1.993333333333331) node[anchor=north west] {$\mathbf{c_1}(-2,4,-2,2)$};
\draw (6.908888888888876,3.175555555555553) node[anchor=north west] {$\mathbf{b_1}(2,1,0,1)$};
\draw (6.482222222222211,4.348888888888886) node[anchor=north west] {$\mathbf{d_1}(6,-2,2,0)$};
\draw (-2.1,0.33) node[anchor=north west] {$\mathbf{o} (0,0,0,0)$};
\draw (1.0511111111111064,4.766666666666664) node[anchor=north west] {$\mathbf{a_1}$};
\draw (0.66,5.7) node[anchor=north west] {$\mathbf{a_2}$};
\draw (-0.36222222222222533,6.091111111111108) node[anchor=north west] {$\mathbf{a_3}$};
\draw (-1.4555555555555575,5.682222222222219) node[anchor=north west] {$\mathbf{a_4}$};
\draw (-1.8,3.9044444444444415) node[anchor=north west] {$-\mathbf{a_2}$};
\draw (0.5,3.9222222222222194) node[anchor=north west] {$-\mathbf{a_4}$};
\draw (-0.7,3.3) node[anchor=north west] {$-\mathbf{a_3}$};
\draw (-2.2,4.8555555555555525) node[anchor=north west] {$-\mathbf{a_1}$};
\draw (-0.5,4.5) node[anchor=north west] {$\mathbf{o}$};
\draw [line width=1.pt,dash pattern=on 2pt off 2pt] (0.,0.)-- (5.656854249492381,5.656854249492381);
\draw [line width=1.pt,dash pattern=on 1pt off 1pt] (0.,0.)-- (6.343145750507619,4.);
\draw [line width=1.pt,dash pattern=on 1pt off 1pt] (0.,0.)-- (6.82842712474619,2.8284271247461903);
\draw [line width=1.pt,dash pattern=on 1pt off 1pt] (0.,0.)-- (7.313708498984761,1.6568542494923806);
\draw [line width=1.pt,dash pattern=on 2pt off 2pt] (0.,0.)-- (8.,0.);
\draw [line width=1.pt] (5.656854249492381,5.656854249492381)-- (8.,0.);
\draw [line width=1.pt] (0.98,4.455555555555551)-- (0.6545630902072708,5.2412297568739366);
\draw [line width=1.pt] (0.6545630902072708,5.2412297568739366)-- (-0.13111111111111437,5.566666666666661);
\draw [line width=1.pt] (-0.13111111111111437,5.566666666666661)-- (-0.9167853124294996,5.2412297568739366);
\draw [line width=1.pt] (-0.9167853124294996,5.2412297568739366)-- (-1.2422222222222243,4.455555555555551);
\draw [line width=1.pt] (-1.2422222222222243,4.455555555555551)-- (-0.9167853124294998,3.669881354237166);
\draw [line width=1.pt] (-0.9167853124294998,3.669881354237166)-- (-0.1311111111111146,3.3444444444444414);
\draw [line width=1.pt] (-0.1311111111111146,3.3444444444444414)-- (0.6545630902072705,3.669881354237166);
\draw [line width=1.pt] (0.6545630902072705,3.669881354237166)-- (0.98,4.455555555555551);
\draw [line width=1.pt,dash pattern=on 2pt off 2pt] (-0.13111111111111448,4.455555555555551)-- (0.98,4.455555555555551);
\draw [line width=1.pt,dash pattern=on 2pt off 2pt] (-0.13111111111111448,4.455555555555551)-- (0.6545630902072708,5.2412297568739366);
\begin{scriptsize}
\draw [fill=ffffff] (0.,0.) circle (3.0pt);
\draw [fill=ffqqqq] (8.,0.) circle (3.0pt);
\draw [fill=ffqqqq] (5.656854249492381,5.656854249492381) circle (3.0pt);
\draw [fill=ududff] (7.313708498984761,1.6568542494923806) circle (2.5pt);
\draw [fill=ududff] (6.82842712474619,2.8284271247461903) circle (2.5pt);
\draw [fill=ududff] (6.343145750507619,4.) circle (2.5pt);
\draw [fill=ffffff] (-0.13111111111111448,4.455555555555551) circle (2.5pt);
\draw [fill=ffqqqq] (0.98,4.455555555555551) circle (2.5pt);
\draw [fill=ffqqqq] (0.6545630902072708,5.2412297568739366) circle (2.5pt);
\draw [fill=ffqqqq] (-0.13111111111111437,5.566666666666661) circle (2.5pt);
\draw [fill=ffqqqq] (-0.9167853124294996,5.2412297568739366) circle (2.5pt);
\draw [fill=ffqqqq] (-1.2422222222222243,4.455555555555551) circle (2.5pt);
\draw [fill=ffqqqq] (-0.9167853124294998,3.669881354237166) circle (2.5pt);
\draw [fill=ffqqqq] (-0.1311111111111146,3.3444444444444414) circle (2.5pt);
\draw [fill=ffqqqq] (0.6545630902072705,3.669881354237166) circle (2.5pt);
\end{scriptsize}
\end{tikzpicture}
\caption{ Defining points $\mathbf{b}_1, \mathbf{c}_1$, and $\mathbf{d}_1$}
\label{fig1}
\end{figure}

In order to do this, consider the points $\mathbf{b}_1$, $\mathbf{c}_1$ and $\mathbf{d}_1$ which divide the side $\mathbf{a}_1\mathbf{a}_2$ in the ratios $1/2$, $1-1/\sqrt{2}$, and $1/\sqrt{2}$, respectively. See Figure \ref{fig1}.

It can be easily checked that
\begin{equation*}
\mathbf{b}_1=(2,1,0,1), \mathbf{c}_1=(-2,4,-2,2), \mathbf{d}_1=(6,-2,2,0).
\end{equation*}


Using equations \eqref{octagonorbit1} and \eqref{octagonorbit2} one can readily compute the octagonal orbits of these three points. We list them below for easy future reference.

\begin{align*}
\langle \mathbf{a}_1\rangle = \{ & (4, 0, 0, 0), (0, 2, 0, 2), (0, 0, 4, 0), (0, -2, 0, 2), \\
&(-4, 0, 0, 0), (0, -2, 0, -2), (0, 0, -4, 0), (0, 2, 0, -2)\},\\
\langle \mathbf{b}_1\rangle=\{ &(2, 1, 0, 1), (0, 1, 2, 1), (0, -1, 2, 1), (-2, -1, 0, 1), \\
&(-2, -1, 0, -1), (0, -1, -2, -1), (0, 1, -2, -1), (2, 1, 0, -1)\},
\end{align*}
\begin{align*}
\langle \mathbf{c}_1\rangle=\{ &(-2, 4, -2, 2), (2, 0, 6, -2), (2, -2, -2, 4), (-6, 2, 2, 0),\\
&(2, -4, 2, -2), (-2, 0, -6, 2), (-2, 2, 2, -4), (6, -2, -2, 0)\},\\
\langle \mathbf{d}_1\rangle=\{ &(6, -2, 2, 0), (-2, 2, -2, 4), (-2, 0, 6, -2), (2, -4, -2, 2),\\
& (-6, 2, -2, 0), (2, -2, 2, -4), (2, 0, -6, 2), (-2, 4, 2, -2)\}.
\end{align*}
Since all points in $\langle \mathbf{a}_1 \rangle$, $\langle \mathbf{b}_1 \rangle$, $\langle \mathbf{c}_1 \rangle$, and
$\langle \mathbf{d}_1 \rangle$ lie on the boundary of the regular octagon, they can be regarded as unit vectors in the octagon metric.


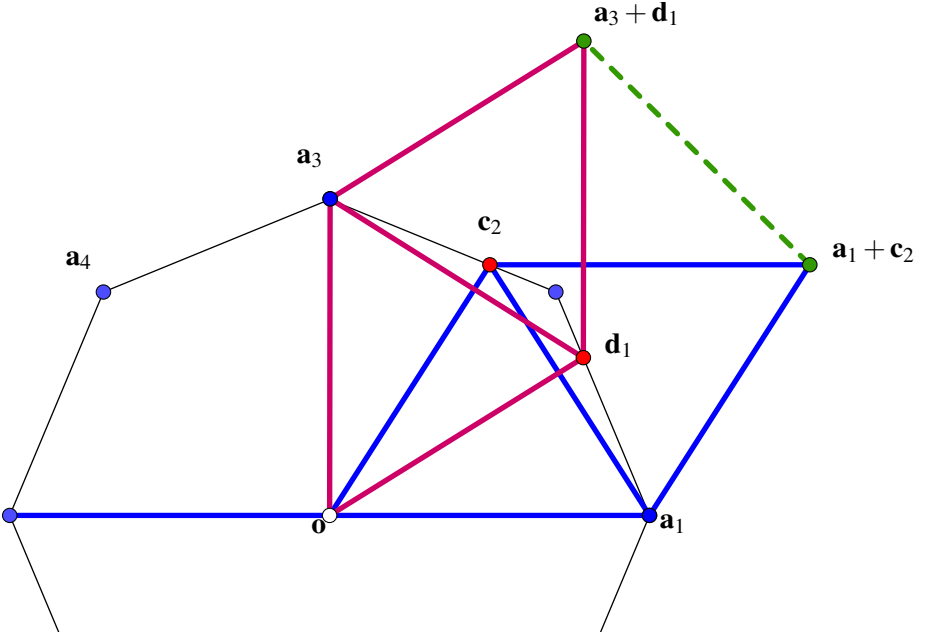
\begin{figure}[htp]
\centering
\begin{tikzpicture}[line cap=round,line join=round,x=1.0cm,y=1.0cm,scale=0.6]
\clip(2.76,-5.58) rectangle (23.360000000000134,8.36);
\draw [line width=0.5pt] (14.949747468305834,1.9497474683058318)-- (10.,4.);
\draw [line width=0.5pt] (10.,4.)-- (5.050252531694169,1.9497474683058336);
\draw [line width=0.5pt] (5.050252531694169,1.9497474683058336)-- (3.,-3.);
\draw [line width=0.5pt] (3.,-3.)-- (5.0502525316941655,-7.94974746830583);
\draw [line width=0.5pt] (5.0502525316941655,-7.94974746830583)-- (10.,-10.);
\draw [line width=0.5pt] (10.,-10.)-- (14.94974746830583,-7.949747468305835);
\draw [line width=0.5pt] (14.94974746830583,-7.949747468305835)-- (17.,-3.);
\draw [line width=0.5pt] (17.,-3.)-- (14.949747468305834,1.9497474683058318);
\draw [line width=2pt,color=qqqqff] (10.,-3.)-- (17.,-3.);
\draw [line width=2pt,color=qqqqff] (10.,-3.)-- (13.507588386092959,2.5533957441303943);
\draw [line width=2pt,color=qqqqff] (10.,-3.)-- (3.,-3.);
\draw [line width=2pt,color=qqqqff] (17.,-3.)-- (13.507588386092959,2.5533957441303943);
\draw [line width=2pt,color=qqqqff] (13.507588386092959,2.5533957441303943)-- (20.50758838609296,2.5533957441303947);
\draw [line width=2pt,color=qqqqff] (17.,-3.)-- (20.50758838609296,2.5533957441303947);
\draw [line width=2pt,color=ccqqww] (10.010710678118656,4.015533905932738)-- (10.,-3.);
\draw [line width=2pt,color=ccqqww] (10.,-3.)-- (15.55239069349651,0.4948380207782339);
\draw [line width=2pt,color=ccqqww] (15.55239069349651,0.4948380207782339)-- (10.010710678118656,4.015533905932738);
\draw [line width=2pt,color=ccqqww] (10.010710678118656,4.015533905932738)-- (15.563101371615168,7.510371926710972);
\draw [line width=2pt,color=ccqqww] (15.563101371615168,7.510371926710972)-- (15.55239069349651,0.4948380207782339);
\draw [line width=2pt,dash pattern=on 5pt off 5pt,color=ttzzqq] (20.50758838609296,2.5533957441303947)-- (15.563101371615168,7.510371926710972);
\draw (17.0000000000001,-2.76) node[anchor=north west] {$\mathbf{a}_1$};
\draw (9.06,5.34) node[anchor=north west] {$\mathbf{a}_3$};
\draw (9.38,-2.86) node[anchor=north west] {$\mathbf{o}$};
\draw (15.8,1.26) node[anchor=north west] {$\mathbf{d}_1$};
\draw (15.56,8.62) node[anchor=north west] {$\mathbf{a}_3+\mathbf{d}_1$};
\draw (13.02,3.86) node[anchor=north west] {$\mathbf{c}_2$};
\draw (20.78000000000012,3.32) node[anchor=north west] {$\mathbf{a}_1+\mathbf{c}_2$};
\draw (4.,3.08) node[anchor=north west] {$\mathbf{a}_4$};
\begin{scriptsize}
\draw [fill=ffffff] (10.,-3.) circle (4.5pt);
\draw [fill=ududff] (17.,-3.) circle (4.5pt);
\draw [fill=ududff] (14.949747468305834,1.9497474683058318) circle (4.5pt);
\draw [fill=ududff] (5.050252531694169,1.9497474683058336) circle (4.5pt);
\draw [fill=ududff] (3.,-3.) circle (4.5pt);
\draw [fill=ududff] (5.0502525316941655,-7.94974746830583) circle (4.5pt);
\draw [fill=ududff] (10.,-10.) circle (4.5pt);
\draw [fill=ududff] (14.94974746830583,-7.949747468305835) circle (4.5pt);
\draw [fill=qqqqff] (17.,-3.) circle (4.5pt);
\draw [fill=ffqqqq] (13.507588386092959,2.5533957441303943) circle (4.5pt);
\draw [fill=ttzzqq] (20.50758838609296,2.5533957441303947) circle (4.5pt);
\draw [fill=ffqqqq] (15.55239069349651,0.4948380207782339) circle (4.5pt);
\draw [fill=qqqqff] (10.010710678118656,4.015533905932738) circle (4.5pt);
\draw [fill=qqqqff] (10.010710678118656,4.015533905932738) circle (4.5pt);
\draw [fill=ttzzqq] (15.563101371615168,7.510371926710972) circle (4.5pt);
\end{scriptsize}
\end{tikzpicture}
\caption{ A Moser spindle in the regular octagon metric}
\label{fig2}
\end{figure}
Consider the points $\mathbf{o}=(0,0,0,0)$, $\mathbf{a}_1=(4,0,0,0)$, $\mathbf{a}_3=(0,0,4,0)$, $\mathbf{c}_2=(2,0,6,-2)$, $\mathbf{d}_1=(6,-2,2,0)$, $\mathbf{a}_1+\mathbf{c}_2=(6,0,6,-2)$, and $\mathbf{a}_3+\mathbf{d}_1=(6,-2,6,0)$.

\noindent One can check that
$\mathbf{c}_2-\mathbf{a}_1=(-2,0,6,-2)=\mathbf{d}_3$, $\mathbf{a}_3-\mathbf{d}_1=(-6,2,2,0)=\mathbf{c}_4$, and
$(\mathbf{a}_3+\mathbf{d}_1)-(\mathbf{a}_1+\mathbf{c}_2)=(0,-2,0,2)=\mathbf{a}_4$. It follows that all bold segments in figure \ref{fig2} are unit vectors in the regular octagon metric. Hence, the Moser spindle can be embedded as a unit-distance graph in the regular octagon metric; it follows that the chromatic number of this plane is at least 4.

We next construct a $5$-chromatic unit distance graph.

Let $U_{32} = \langle \mathbf{a}_1\rangle \cup \langle \mathbf{b}_1\rangle \cup \langle \mathbf{c}_1\rangle \cup\langle \mathbf{d}_1\rangle$ be the set of generating vectors, and consider the $2$-fold Minkowski sum
\begin{equation}
U_{32}+U_{32} =\{ \mathbf{x}+\mathbf{y} \,\, | \,\, \mathbf{x}\in U_{32}, \mathbf{y} \in U_{32}\}.
\end{equation}

A straightforward computation shows that the unit-distance graph with vertex set $U_{32}+U_{32}$ has $465$ vertices and $2368$ edges in $U_{32}$.
The chromatic number of this graph is still $4$.

There are however many pairs of vertices $\mathbf{u}, \mathbf{v} \in U_{32}+U_{32}$ which are unit-distance apart in the regular octagon metric despite the fact that $\mathbf{u}-\mathbf{v}\notin U_{32}$. Indeed, consider the following vectors
\begin{align*}
\mathbf{e}_1:&=(\mathbf{a}_1+\mathbf{a}_3)-(\mathbf{c}_1+\mathbf{c}_4)=(12,-6,4,-2),\\
\mathbf{f}_1:&=(\mathbf{c}_1+\mathbf{d}_2)-(\mathbf{a}_1+\mathbf{a}_4)=(-8,8,-4,4),\\
\mathbf{g}_1:&=(\mathbf{a}_2+\mathbf{c}_1)-(\mathbf{a}_1+\mathbf{b}_3)=(-6,7,-4,3),\\
\mathbf{h}_1:&=(\mathbf{a}_3+\mathbf{b}_1)-(\mathbf{c}_1+\mathbf{c}_4)=(10,-5,4,-1),\\
\mathbf{i}_1:&=(\mathbf{a}_1+\mathbf{a}_3)-(\mathbf{b}_1+\mathbf{c}_4)=(8,-3,2,-1),\\
\mathbf{j}_1:&=(\mathbf{a}_2+\mathbf{d}_2)-(\mathbf{a}_4+\mathbf{b}_1)=(-4,5,-2,3).\\
\end{align*}

\begin{figure}[htp]
\centering
\begin{tikzpicture}[line cap=round,line join=round,x=1.0cm,y=1.0cm,scale=1]
\clip(-2.0066666666666677,-0.3) rectangle (11,6.588888888888888);
\draw (5.726666666666657,6.224444444444443) node[anchor=north west] {$\mathbf{a_2}(0,2,0,2)$};
\draw (8.073333333333322,0.3666666666666666) node[anchor=north west] {$\mathbf{a_1} (4,0,0,0)$};
\draw (7.113333333333322,2.7488888888888883) node[anchor=north west] {$\mathbf{e_1}(12,-6,4,-2)$};
\draw (6.713333333333322,3.7088888888888882) node[anchor=north west] {$\mathbf{f_1}(-8,8,-4,4)$};
\draw [line width=1.pt,dash pattern=on 2pt off 2pt] (0.,0.)-- (5.656854249492381,5.656854249492381);
\draw [line width=1.pt,dash pattern=on 2pt off 2pt] (0.,0.)-- (8.,0.);
\draw [line width=1.pt] (5.656854249492381,5.656854249492381)-- (8.,0.);
\draw (-1.8644444444444455,0.3577777777777777) node[anchor=north west] {$\mathbf{o} (0,0,0,0)$};
\draw [line width=1.pt] (0.98,4.455555555555551)-- (0.6545630902072708,5.2412297568739366);
\draw [line width=1.pt] (0.6545630902072708,5.2412297568739366)-- (-0.13111111111111437,5.566666666666661);
\draw [line width=1.pt] (-0.13111111111111437,5.566666666666661)-- (-0.9167853124294996,5.2412297568739366);
\draw [line width=1.pt] (-0.9167853124294996,5.2412297568739366)-- (-1.2422222222222243,4.455555555555551);
\draw [line width=1.pt] (-1.2422222222222243,4.455555555555551)-- (-0.9167853124294998,3.669881354237166);
\draw [line width=1.pt] (-0.9167853124294998,3.669881354237166)-- (-0.1311111111111146,3.3444444444444414);
\draw [line width=1.pt] (-0.1311111111111146,3.3444444444444414)-- (0.6545630902072705,3.669881354237166);
\draw [line width=1.pt] (0.6545630902072705,3.669881354237166)-- (0.98,4.455555555555551);
\draw (1.0511111111111067,4.766666666666666) node[anchor=north west] {$\mathbf{a_1}$};
\draw (0.66,5.7) node[anchor=north west] {$\mathbf{a_2}$};
\draw (-0.3622222222222251,6.09111111111111) node[anchor=north west] {$\mathbf{a_3}$};
\draw (-1.455555555555557,5.682222222222221) node[anchor=north west] {$\mathbf{a_4}$};
\draw (-1.66,3.9044444444444437) node[anchor=north west] {$-\mathbf{a_2}$};
\draw (0.5888888888888849,3.9222222222222216) node[anchor=north west] {$-\mathbf{a_4}$};
\draw (-0.5666666666666693,3.4511111111111106) node[anchor=north west] {$-\mathbf{a_3}$};
\draw (-1.9977777777777788,4.855555555555554) node[anchor=north west] {$-\mathbf{a_1}$};
\draw (-0.19333333333333638,4.58) node[anchor=north west] {$\mathbf{o}$};
\draw [line width=1.pt,dash pattern=on 2pt off 2pt] (-0.13111111111111448,4.455555555555551)-- (0.98,4.455555555555551);
\draw [line width=1.pt,dash pattern=on 2pt off 2pt] (-0.13111111111111448,4.455555555555551)-- (0.6545630902072708,5.2412297568739366);
\draw [line width=1.pt] (0.,0.)-- (7.029437251522857,2.3431457505076194);
\draw [line width=1.pt] (0.,0.)-- (6.627416997969522,3.313708498984761);
\draw (7.886666666166654,0.9177777777777776) node[anchor=north west] {$\mathbf{g}_1(-6,7,-4,3)$};
\draw [line width=1.pt] (0.,0.)-- (7.798989873223331,0.4852813742385713);
\draw (5.913333333333323,5.584444444444443) node[anchor=north west] {$\mathbf{h}_1(10, -5, 4, -1)$};
\draw [line width=1.pt] (0.,0.)-- (5.857864376269049,5.17157287525381);
\draw (7.628888888888877,1.584444444444444) node[anchor=north west] {$\mathbf{i}_1(8,-3,2,-1)$};
\draw (6.26,4.873333333333332) node[anchor=north west] {$\mathbf{j}_1(-4,5,-2,3)$};
\draw [line width=1.pt] (0.,0.)-- (6.142135623730951,4.485281374238571);
\draw [line width=1.pt] (0.,0.)-- (7.514718625761429,1.1715728752538097);
\begin{scriptsize}
\draw [fill=ffffff] (0.,0.) circle (3.0pt);
\draw [fill=ffqqqq] (8.,0.) circle (3.0pt);
\draw [fill=ffqqqq] (5.656854249492381,5.656854249492381) circle (3.0pt);
\draw [fill=ffffff] (-0.13111111111111448,4.455555555555551) circle (2.5pt);
\draw [fill=ffqqqq] (0.98,4.455555555555551) circle (2.5pt);
\draw [fill=ffqqqq] (0.6545630902072708,5.2412297568739366) circle (2.5pt);
\draw [fill=ffqqqq] (-0.13111111111111437,5.566666666666661) circle (2.5pt);
\draw [fill=ffqqqq] (-0.9167853124294996,5.2412297568739366) circle (2.5pt);
\draw [fill=ffqqqq] (-1.2422222222222243,4.455555555555551) circle (2.5pt);
\draw [fill=ffqqqq] (-0.9167853124294998,3.669881354237166) circle (2.5pt);
\draw [fill=ffqqqq] (-0.1311111111111146,3.3444444444444414) circle (2.5pt);
\draw [fill=ffqqqq] (0.6545630902072705,3.669881354237166) circle (2.5pt);
\draw [fill=ududff] (7.029437251522857,2.3431457505076194) circle (2.5pt);
\draw [fill=ududff] (6.627416997969522,3.313708498984761) circle (2.5pt);
\draw [fill=ududff] (7.798989873223331,0.4852813742385713) circle (2.5pt);
\draw [fill=ududff] (5.857864376269049,5.17157287525381) circle (2.5pt);
\draw [fill=ududff] (7.514718625761429,1.1715728752538097) circle (2.5pt);
\draw [fill=ududff] (6.142135623730951,4.485281374238571) circle (2.5pt);
\end{scriptsize}
\end{tikzpicture}
\caption{ Accidental unit vectors in the regular octagon metric.}
\label{fig3}
\end{figure}
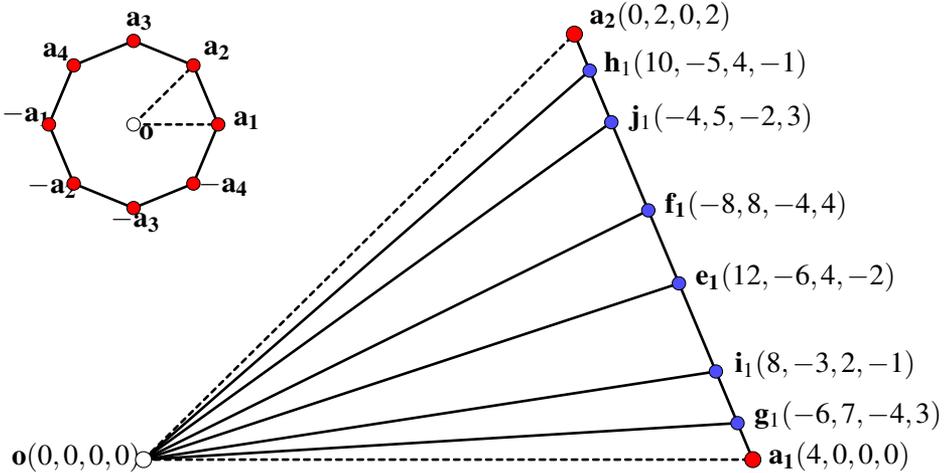

While none of these six vectors belongs to $U_{32}$, they are all unit vectors in the octagon metric.
Indeed, it can be verified without difficulty that $\mathbf{e}_1, \mathbf{f}_1,\mathbf{g}_1, \mathbf{h}_1,\mathbf{i}_1, \mathbf{j}_1$ divide the edge $\mathbf{a}_1\mathbf{a}_2$ of the octagon in the ratios $\sqrt{2}-1, 2-\sqrt{2}, 3/2-\sqrt{2},\\ \sqrt{2}-1/2, (\sqrt{2}-1)/2$, and $(3-\sqrt{2})/2$, respectively. See figure \ref{fig3}.

\newpage
Using equations \eqref{octagonorbit1} and \eqref{octagonorbit2} one can easily compute the octagonal orbits of the accidental unit vectors $\mathbf{e}_1, \mathbf{f}_1,\mathbf{g}_1, \mathbf{h}_1,\mathbf{i}_1, \mathbf{j}_1$; due to space constraints we omit listing them here.
Let $V_{48}$ be the union of these orbits
\begin{equation*}
V_{48}:=\langle \mathbf{e}_1\rangle \cup \langle \mathbf{f}_1\rangle \cup \langle \mathbf{g}_1\rangle \cup\langle \mathbf{h}_1\rangle \cup \langle \mathbf{i}_1\rangle \cup\langle \mathbf{j}_1\rangle.
\end{equation*}
The set $V_{48}$ is the set of accidental vectors we mentioned in the introduction.
We now have a set of $80$ unit vectors: $W_{80}:= U_{32} \cup V_{48}$. Returning to the graph with vertex set $U_{32}+U_{32}$, besides the
$2368$ unit edges in $U_{32}$ there are another $1072$ unit edges in $V_{48}$. This graph has chromatic number $5$.

It may be desirable to reduce the order of the $5$-chromatic graph. Our best construction is presented in the following

\begin{thm}
Let $S$ consist of the following $15$ points in the Minkowski plane endowed with the regular octagon metric:
\begin{align*}
S:=\{&(-4, 4, 0, 0), (-4, 5, 2, -1), (-2, 2, 6, -4), (-2, 3, -4, 3), (-2, 4, -2, 2),\\
 &(0, 1, 2, -1),(0, 2, 4, -2), (2, -1, -4, 3), (2, 0, -2, 2), (2, 1, 0, 1),\\
 &(2, 2, 2, 0),(4, 0, -4, 4), (4, 0, 0, 0), (6, -2, 2, 0), (6, -1, 0, 1)\}
\end{align*}
Then the graph $G_{120}$, whose vertex set is $\bigcup_{i=1}^{15} \langle S_i \rangle$, has $120$ vertices, $704$ edges in $W_{80}$ and has chromatic number $5$. Hence, $\chi(\mathbf{R}^2, C)\ge 5$ when $C$ is a regular octagon.
\end{thm}
\begin{proof}
For verification purposes we provide the edge distribution over the unit vector orbits $\langle\mathbf{a}_1\rangle, \langle\mathbf{b}_1\rangle, \langle\mathbf{c}_1\rangle, \langle\mathbf{d}_1\rangle, \langle\mathbf{e}_1\rangle, \langle\mathbf{f}_1\rangle, \langle\mathbf{g}_1\rangle, \langle\mathbf{h}_1\rangle,\langle\mathbf{i}_1\rangle, \langle\mathbf{j}_1\rangle$:
\begin{equation*}
160, 128, 136, 128, 16, 8, 32, 8, 32, 56.
\end{equation*}
\end{proof}

The chromatic number verification for all of the graphs described in this paper
can be completed in under a minute using Maple.
It can be checked that the graph contains $24$ Moser spindles.
It is depicted in the two figures below.
Figure~\ref{octexpected} shows the edges of $G_{120}$ corresponding to vectors in $U_{32}$,
while Figure~\ref{octunexpected} shows the edges of $G_{120}$ corresponding to vectors
in $V_{48}$.

\begin{figure}[htp]
\centering
\begin{tikzpicture}[line width=1pt,scale=0.8]
\tikzstyle{every node}=[draw=black,fill=yellow!50!white,thick,
  shape=circle,minimum height=0.2cm,inner sep=1];
\tikzstyle{every path}=[draw=white!10!black];

\input{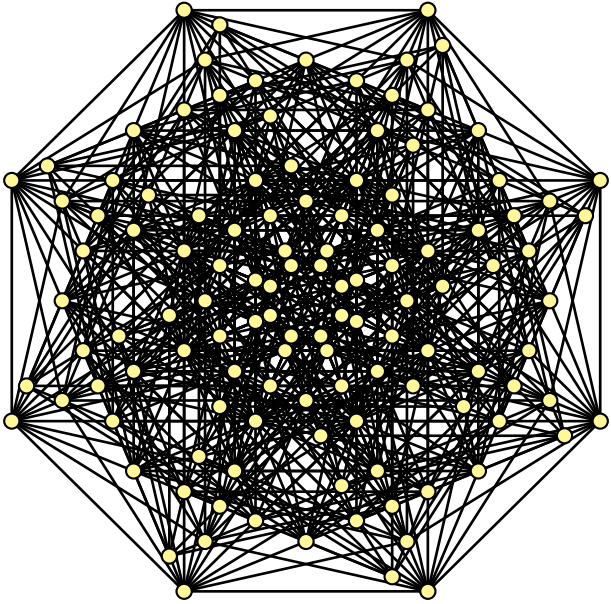}

\end{tikzpicture}
\caption{The expected edges in the $5$-chromatic graph $G_{120}$ for the regular octagon metric}
\label{octexpected}
\end{figure}

\begin{figure}[htp]
\centering
\begin{tikzpicture}[line width=1pt,scale=0.8]
\tikzstyle{every node}=[draw=black,fill=yellow!50!white,thick,
  shape=circle,minimum height=0.2cm,inner sep=1];
\tikzstyle{every path}=[draw=white!10!black];

\input{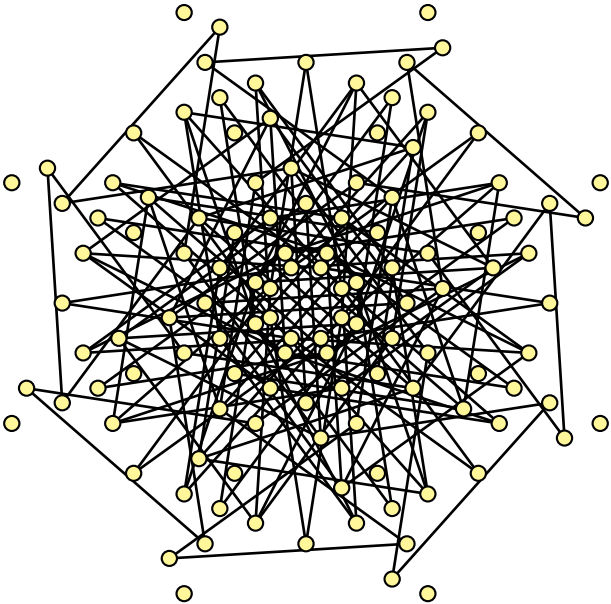}

\end{tikzpicture}
\caption{The accidental edges in the $5$-chromatic graph  $G_{120}$ for the regular octagon metric}
\label{octunexpected}
\end{figure}

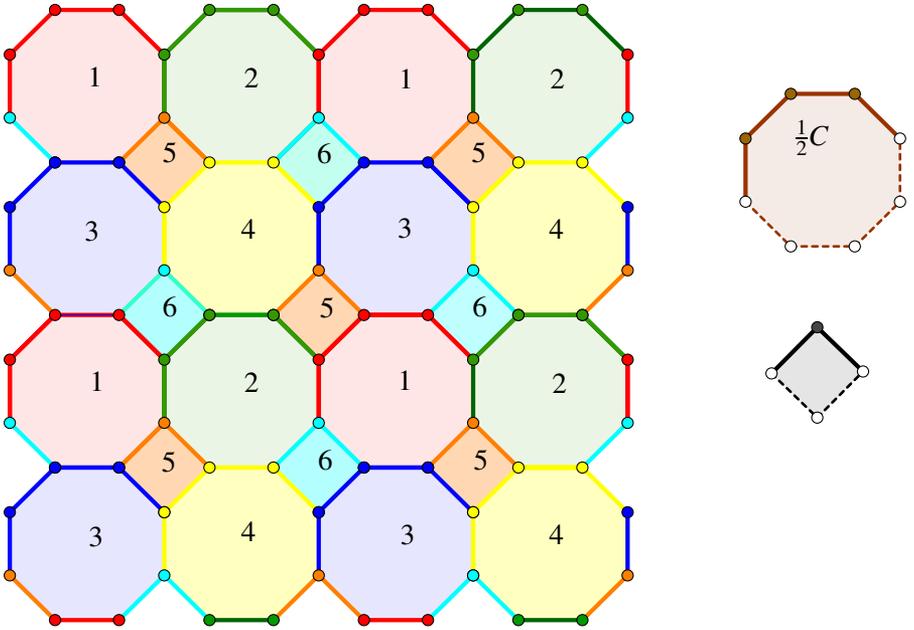
\begin{figure}[htp]
\centering
\begin{tikzpicture}[line cap=round,line join=round,x=1.0cm,y=1.0cm,scale=0.42]
\clip(11.148333333333358,-10.08) rectangle (39.738333333333394,9.9);
\fill[line width=0.pt,color=ffqqqq,fill=ffqqqq,fill opacity=0.10000000149011612] (16.414213562373057,6.071067811865503) -- (16.414213562373043,8.071067811865507) -- (15.,9.485281374238548) -- (13.,9.485281374238548) -- (11.585786437626936,8.071067811865476) -- (11.585786437626934,6.071067811865506) -- (13.,4.656854249492376) -- (15.,4.6568542494923975) -- cycle;
\fill[line width=0.pt,color=ffqqqq,fill=ffqqqq,fill opacity=0.10000000149011612] (21.242640687119156,6.071067811865542) -- (21.242640687119156,8.07106781186551) -- (22.6568542494922,9.485281374238566) -- (24.656854249492163,9.485281374238589) -- (26.071067811865248,8.071067811865536) -- (26.071067811865273,6.071067811865571) -- (24.65685424949222,4.6568542494924845) -- (22.656854249492255,4.656854249492461) -- cycle;
\fill[line width=0.0pt,color=ffqqqq,fill=ffqqqq,fill opacity=0.10000000149011612] (13.,-5.) -- (15.,-5.) -- (16.414213562373096,-3.585786437626905) -- (16.414213562373096,-1.5857864376269026) -- (15.,-0.17157287525380602) -- (13.,-0.1715728752538057) -- (11.5857864376269,-1.5857864376269015) -- (11.585786437626899,-3.585786437626904) -- cycle;
\fill[line width=0.0pt,color=ffqqqq,fill=ffqqqq,fill opacity=0.10000000149011612] (26.071067811865518,-1.5857864376268749) -- (24.656854249492298,-0.17157287525375875) -- (22.656854249492305,-0.17157287525388165) -- (21.242640687119305,-1.5857864376270585) -- (21.24264068711943,-3.585786437627049) -- (22.656854249492426,-5.) -- (24.656854249492444,-5.) -- (26.071067811865536,-3.585786437626891) -- cycle;
\fill[line width=0.0pt,color=ttzzqq,fill=ttzzqq,fill opacity=0.10000000149011612] (16.414213562373057,6.071067811865503) -- (17.828427124746117,4.656854249492467) -- (19.828427124746085,4.6568542494924685) -- (21.242640687119156,6.071067811865542) -- (21.242640687119156,8.07106781186551) -- (19.82842712474608,9.485281374238582) -- (17.828427124746113,9.48528137423858) -- (16.414213562373043,8.071067811865507) -- cycle;
\fill[line width=0.pt,color=ttzzqq,fill=ttzzqq,fill opacity=0.10000000149011612] (30.89949493661139,8.071067811865671) -- (29.485281374238248,9.485281374238703) -- (27.485281374238262,9.485281374238625) -- (26.071067811865248,8.071067811865536) -- (26.071067811865273,6.071067811865571) -- (27.485281374238593,4.65685424949263) -- (29.485281374238593,4.6568542494926035) -- (30.899494936611468,6.0710678118656825) -- cycle;
\fill[line width=0.0pt,color=ttzzqq,fill=ttzzqq,fill opacity=0.10000000149011612] (21.24264068711943,-3.585786437627049) -- (21.242640687119305,-1.5857864376270585) -- (19.8284271247462,-0.17157287525380616) -- (17.828427124746195,-0.17157287525380438) -- (16.414213562373096,-1.5857864376269026) -- (16.414213562373096,-3.585786437626905) -- (17.82842712474619,-5.) -- (19.828427124746195,-5.) -- cycle;
\fill[line width=0.0pt,color=ttzzqq,fill=ttzzqq,fill opacity=0.10000000149011612] (30.899494936611788,-1.5857864376267878) -- (29.485281374238532,-0.17157287525359122) -- (27.48528137423853,-0.17157287525356602) -- (26.071067811865518,-1.5857864376268749) -- (26.071067811865536,-3.585786437626891) -- (27.485281374238706,-5.) -- (29.485281374238745,-5.) -- (30.899494936611834,-3.585786437626825) -- cycle;
\fill[line width=0.0pt,color=qqqqff,fill=qqqqff,fill opacity=0.10000000149011612] (16.414213562373106,1.2426406871193165) -- (16.414213562373085,3.2426406871193167) -- (15.,4.6568542494923975) -- (13.,4.656854249492376) -- (11.585786437626894,3.2426406871192666) -- (11.585786437626915,1.2426406871192657) -- (13.,-0.1715728752538057) -- (15.,-0.17157287525380602) -- cycle;
\fill[line width=0.0pt,color=ffffqq,fill=ffffqq,fill opacity=0.25] (21.242640687119287,1.2426406871193783) -- (21.24264068711926,3.242640687119373) -- (19.828427124746085,4.6568542494924685) -- (17.828427124746117,4.656854249492467) -- (16.414213562373085,3.2426406871193167) -- (16.414213562373106,1.2426406871193165) -- (17.828427124746195,-0.17157287525380438) -- (19.8284271247462,-0.17157287525380616) -- cycle;
\fill[line width=0.0pt,color=qqqqff,fill=qqqqff,fill opacity=0.10000000149011612] (22.656854249492305,-0.17157287525388165) -- (24.656854249492298,-0.17157287525375875) -- (26.071067811865483,1.242640687119466) -- (26.071067811865436,3.2426406871194624) -- (24.65685424949222,4.6568542494924845) -- (22.656854249492255,4.656854249492461) -- (21.24264068711926,3.242640687119373) -- (21.242640687119287,1.2426406871193783) -- cycle;
\fill[line width=0.pt,color=ffffqq,fill=ffffqq,fill opacity=0.25] (30.89949493661167,3.24264068711949) -- (29.485281374238593,4.6568542494926035) -- (27.485281374238593,4.65685424949263) -- (26.071067811865436,3.2426406871194624) -- (26.071067811865483,1.242640687119466) -- (27.48528137423853,-0.17157287525356602) -- (29.485281374238532,-0.17157287525359122) -- (30.899494936611646,1.2426406871194882) -- cycle;
\fill[line width=0.pt,color=qqqqff,fill=qqqqff,fill opacity=0.10000000149011612] (11.585786437626902,-8.414213562373094) -- (13.,-9.82842712474619) -- (15.,-9.828427124746192) -- (16.414213562373092,-8.414213562373098) -- (16.414213562373092,-6.4142135623730985) -- (15.,-5.) -- (13.,-5.) -- (11.585786437626904,-6.414213562373095) -- cycle;
\fill[line width=0.pt,color=qqqqff,fill=qqqqff,fill opacity=0.10000000149011612] (26.071067811865532,-6.414213562373101) -- (24.656854249492444,-5.) -- (22.656854249492426,-5.) -- (21.242640687119298,-6.414213562373112) -- (21.242640687119298,-8.414213562373124) -- (22.65685424949246,-9.828427124746172) -- (24.656854249492447,-9.828427124746172) -- (26.071067811865532,-8.414213562373087) -- cycle;
\fill[line width=0.pt,color=ffffqq,fill=ffffqq,fill opacity=0.25] (17.82842712474618,-9.828427124746225) -- (19.82842712474619,-9.828427124746225) -- (21.242640687119298,-8.414213562373124) -- (21.242640687119298,-6.414213562373112) -- (19.828427124746195,-5.) -- (17.82842712474619,-5.) -- (16.414213562373092,-6.4142135623730985) -- (16.414213562373092,-8.414213562373098) -- cycle;
\fill[line width=0.pt,color=ffffqq,fill=ffffqq,fill opacity=0.25] (30.899494936611895,-8.414213562373105) -- (30.899494936611887,-6.414213562373092) -- (29.485281374238745,-5.) -- (27.485281374238706,-5.) -- (26.071067811865532,-6.414213562373101) -- (26.071067811865532,-8.414213562373087) -- (27.485281374238784,-9.828427124746222) -- (29.485281374238795,-9.828427124746213) -- cycle;
\fill[line width=0.pt,color=ffxfqq,fill=ffxfqq,fill opacity=0.30000001192092896] (16.414213562373057,6.071067811865503) -- (15.,4.6568542494923975) -- (16.414213562373085,3.2426406871193167) -- (17.828427124746117,4.656854249492467) -- cycle;
\fill[line width=0.pt,color=ffxfqq,fill=ffxfqq,fill opacity=0.30000001192092896] (19.8284271247462,-0.17157287525380616) -- (21.242640687119305,-1.5857864376270585) -- (22.656854249492305,-0.17157287525388165) -- (21.242640687119287,1.2426406871193783) -- cycle;
\fill[line width=0.0pt,color=ffxfqq,fill=ffxfqq,fill opacity=0.30000001192092896] (26.071067811865273,6.071067811865571) -- (24.65685424949222,4.6568542494924845) -- (26.071067811865436,3.2426406871194624) -- (27.485281374238593,4.65685424949263) -- cycle;
\fill[line width=0.pt,color=ffxfqq,fill=ffxfqq,fill opacity=0.30000001192092896] (15.,-5.) -- (16.414213562373092,-6.4142135623730985) -- (17.82842712474619,-5.) -- (16.414213562373096,-3.585786437626905) -- cycle;
\fill[line width=0.pt,color=ffxfqq,fill=ffxfqq,fill opacity=0.30000001192092896] (24.656854249492444,-5.) -- (26.071067811865532,-6.414213562373101) -- (27.485281374238706,-5.) -- (26.071067811865536,-3.585786437626891) -- cycle;
\fill[line width=0.0pt,color=ttffcc,fill=ttffcc,fill opacity=0.30000001192092896] (21.242640687119156,6.071067811865542) -- (19.828427124746085,4.6568542494924685) -- (21.24264068711926,3.242640687119373) -- (22.656854249492255,4.656854249492461) -- cycle;
\fill[line width=0.pt,color=qqffff,fill=qqffff,fill opacity=0.30000001192092896] (15.,-0.17157287525380602) -- (16.414213562373096,-1.5857864376269026) -- (17.828427124746195,-0.17157287525380438) -- (16.414213562373106,1.2426406871193165) -- cycle;
\fill[line width=0.pt,color=qqffff,fill=qqffff,fill opacity=0.25] (26.071067811865483,1.242640687119466) -- (24.656854249492298,-0.17157287525375875) -- (26.071067811865518,-1.5857864376268749) -- (27.48528137423853,-0.17157287525356602) -- cycle;
\fill[line width=0.pt,color=qqffff,fill=qqffff,fill opacity=0.30000001192092896] (19.828427124746195,-5.) -- (21.242640687119298,-6.414213562373112) -- (22.656854249492426,-5.) -- (21.24264068711943,-3.585786437627049) -- cycle;
\fill[line width=0.pt,color=zzttqq,fill=zzttqq,fill opacity=0.10000000149011612] (39.41421356237309,5.414213562373091) -- (38.,6.828427124746184) -- (36.,6.8284271247461845) -- (34.58578643762691,5.41421356237309) -- (34.58578643762691,3.414213562373092) -- (36.,2.) -- (38.,2.) -- (39.41421356237309,3.414213562373095) -- cycle;
\fill[line width=0.pt,color=yqyqyq,fill=yqyqyq,fill opacity=0.20000000298023224] (38.258333333333404,-1.96) -- (36.828333333333404,-0.5596853493024382) -- (35.3983333333334,-2.02) -- (36.828333333333404,-3.420314650697562) -- cycle;
\draw [line width=1.6pt,color=ffqqqq] (11.585786437626934,6.071067811865506)-- (11.585786437626936,8.071067811865476);
\draw [line width=1.6pt,color=ffqqqq] (11.585786437626936,8.071067811865476)-- (13.,9.485281374238548);
\draw [line width=1.6pt,color=ffqqqq] (13.,9.485281374238548)-- (15.,9.485281374238548);
\draw [line width=1.6pt,color=ffqqqq] (15.,9.485281374238548)-- (16.414213562373043,8.071067811865507);
\draw [line width=1.6pt,color=qqqqff] (11.585786437626915,1.2426406871192657)-- (11.585786437626894,3.2426406871192666);
\draw [line width=1.6pt,color=qqqqff] (11.585786437626894,3.2426406871192666)-- (13.,4.656854249492376);
\draw [line width=1.6pt,color=qqqqff] (13.,4.656854249492376)-- (15.,4.6568542494923975);
\draw [line width=1.6pt,color=qqqqff] (15.,4.6568542494923975)-- (16.414213562373085,3.2426406871193167);
\draw [line width=1.6pt,color=qqqqff] (11.585786437626902,-8.414213562373094)-- (11.585786437626904,-6.414213562373095);
\draw [line width=0.0pt,color=qqqqff] (11.585786437626904,-6.414213562373095)-- (13.,-5.);
\draw [line width=1.6pt,color=qqqqff] (13.,-5.)-- (15.,-5.);
\draw [line width=1.6pt,color=qqqqff] (15.,-5.)-- (16.414213562373092,-6.4142135623730985);
\draw [line width=0.0pt,color=qqqqff] (21.24264068711926,3.242640687119373)-- (22.656854249492316,4.656854249492477);
\draw [line width=1.6pt,color=qqqqff] (22.656854249492316,4.656854249492477)-- (24.656854249492312,4.656854249492523);
\draw [line width=1.6pt,color=qqqqff] (21.242640687119287,1.2426406871193783)-- (21.24264068711926,3.242640687119373);
\draw [line width=1.6pt,color=qqqqff] (21.242640687119298,-8.414213562373124)-- (21.242640687119298,-6.414213562373112);
\draw [line width=1.6pt,color=qqqqff] (21.242640687119298,-6.414213562373112)-- (22.656854249492426,-5.);
\draw [line width=1.6pt,color=qqqqff] (22.656854249492426,-5.)-- (24.656854249492444,-5.);
\draw [line width=0.4pt,color=qqqqff] (24.656854249492312,4.656854249492523)-- (26.071067811865436,3.2426406871194624);
\draw [line width=1.6pt,color=qqqqff] (24.656854249492444,-5.)-- (26.071067811865532,-6.414213562373101);
\draw [line width=1.6pt,color=ffqqqq] (21.242640687119156,6.071067811865542)-- (21.242640687119156,8.07106781186551);
\draw [line width=1.6pt,color=ffqqqq] (21.242640687119156,8.07106781186551)-- (22.6568542494922,9.485281374238566);
\draw [line width=1.6pt,color=ffqqqq] (22.6568542494922,9.485281374238566)-- (24.656854249492163,9.485281374238589);
\draw [line width=1.6pt,color=ffqqqq] (24.656854249492163,9.485281374238589)-- (26.071067811865248,8.071067811865536);
\draw [line width=1.6pt,color=ffqqqq] (21.242640687119295,-3.5857864376269077)-- (21.242640687119295,-1.585786437626904);
\draw [line width=1.6pt,color=ffqqqq] (21.242640687119295,-1.585786437626904)-- (22.65685424949238,-0.1715728752537995);
\draw [line width=1.6pt,color=ffqqqq] (22.65685424949238,-0.1715728752537995)-- (24.656854249492397,-0.17157287525378084);
\draw [line width=1.6pt,color=ffqqqq] (24.656854249492397,-0.17157287525378084)-- (26.071067811865518,-1.5857864376268749);
\draw [line width=1.6pt,color=ffqqqq] (11.5857864376269,-3.5857864376269024)-- (11.5857864376269,-1.5857864376269004);
\draw [line width=1.6pt,color=ffqqqq] (11.5857864376269,-1.5857864376269004)-- (13.,-0.17157287525380482);
\draw [line width=0.0pt,color=ffqqqq] (13.,-0.17157287525380482)-- (15.,-0.1715728752538057);
\draw [line width=1.6pt,color=ffqqqq] (15.,-0.1715728752538057)-- (16.414213562373096,-1.5857864376269026);
\draw [line width=1.6pt,color=qqwuqq] (16.414213562373043,8.071067811865507)-- (17.828427124746113,9.48528137423858);
\draw [line width=1.6pt,color=qqwuqq] (17.828427124746113,9.48528137423858)-- (19.82842712474608,9.485281374238582);
\draw [line width=1.6pt,color=qqwuqq] (16.414213562373057,6.071067811865503)-- (16.414213562373043,8.071067811865507);
\draw [line width=1.6pt,color=qqwuqq] (19.82842712474608,9.485281374238582)-- (21.242640687119156,8.07106781186551);
\draw [line width=1.6pt,color=qqwuqq] (16.414213562373096,-1.5857864376269026)-- (17.828427124746195,-0.17157287525380438);
\draw [line width=1.6pt,color=qqwuqq] (16.414213562373096,-3.585786437626905)-- (16.414213562373096,-1.5857864376269026);
\draw [line width=0.0pt,color=qqwuqq] (17.828427124746195,-0.17157287525380438)-- (19.8284271247462,-0.17157287525380616);
\draw [line width=1.6pt,color=qqwuqq] (19.8284271247462,-0.17157287525380616)-- (21.242640687119295,-1.585786437626904);
\draw [line width=1.6pt,color=qqwuqq] (26.071067811865273,6.071067811865571)-- (26.071067811865248,8.071067811865536);
\draw [line width=1.6pt,color=qqwuqq] (26.071067811865248,8.071067811865536)-- (27.485281374238262,9.485281374238625);
\draw [line width=1.6pt,color=qqwuqq] (27.485281374238262,9.485281374238625)-- (29.485281374238248,9.485281374238703);
\draw [line width=1.6pt,color=qqwuqq] (29.485281374238248,9.485281374238703)-- (30.89949493661139,8.071067811865671);
\draw [line width=1.6pt,color=qqwuqq] (26.071067811865536,-3.585786437626891)-- (26.071067811865518,-1.5857864376268749);
\draw [line width=1.6pt,color=qqwuqq] (26.071067811865518,-1.5857864376268749)-- (27.48528137423853,-0.17157287525356602);
\draw [line width=1.6pt,color=qqwuqq] (27.48528137423853,-0.17157287525356602)-- (29.485281374238532,-0.17157287525359122);
\draw [line width=1.6pt,color=qqwuqq] (29.485281374238532,-0.17157287525359122)-- (30.899494936611788,-1.5857864376267878);
\draw [line width=1.6pt,color=ffffqq] (16.414213562373106,1.2426406871193165)-- (16.414213562373085,3.2426406871193167);
\draw [line width=1.6pt,color=ffffqq] (16.414213562373085,3.2426406871193167)-- (17.828427124746156,4.656854249492417);
\draw [line width=1.6pt,color=ffffqq] (17.828427124746156,4.656854249492417)-- (19.82842712474615,4.6568542494924445);
\draw [line width=1.6pt,color=ffffqq] (19.82842712474615,4.6568542494924445)-- (21.24264068711926,3.242640687119373);
\draw [line width=1.6pt,color=ffffqq] (26.071067811865483,1.242640687119466)-- (26.071067811865436,3.2426406871194624);
\draw [line width=1.6pt,color=ffffqq] (26.071067811865436,3.2426406871194624)-- (27.485281374238447,4.656854249492461);
\draw [line width=1.6pt,color=ffffqq] (27.485281374238447,4.656854249492461)-- (29.485281374238436,4.65685424949254);
\draw [line width=1.6pt,color=ffffqq] (29.485281374238436,4.65685424949254)-- (30.89949493661167,3.24264068711949);
\draw [line width=1.6pt,color=ffffqq] (16.414213562373092,-6.4142135623730985)-- (17.82842712474619,-5.);
\draw [line width=1.6pt,color=ffffqq] (16.414213562373092,-8.414213562373098)-- (16.414213562373092,-6.4142135623730985);
\draw [line width=1.6pt,color=ffffqq] (17.82842712474619,-5.)-- (19.828427124746195,-5.);
\draw [line width=1.6pt,color=ffffqq] (19.828427124746195,-5.)-- (21.242640687119298,-6.414213562373112);
\draw [line width=1.6pt,color=ffffqq] (26.071067811865532,-8.414213562373087)-- (26.071067811865532,-6.414213562373101);
\draw [line width=1.6pt,color=ffffqq] (26.071067811865532,-6.414213562373101)-- (27.485281374238706,-5.);
\draw [line width=2.pt,color=ffffqq] (27.485281374238706,-5.)-- (29.485281374238745,-5.);
\draw [line width=1.6pt,color=ffffqq] (29.485281374238745,-5.)-- (30.899494936611887,-6.414213562373092);
\draw [line width=1.6pt,color=qqqqff] (13.,-5.)-- (15.,-5.);
\draw [line width=1.6pt,color=ffxfqq] (15.,-5.)-- (16.414213562373096,-3.585786437626905);
\draw [line width=1.6pt,color=ttzzqq] (16.414213562373096,-3.585786437626905)-- (16.414213562373096,-1.5857864376269026);
\draw [line width=1.6pt,color=ffqqqq] (16.414213562373096,-1.5857864376269026)-- (15.,-0.17157287525380602);
\draw [line width=1.6pt,color=ffqqqq] (15.,-0.17157287525380602)-- (13.,-0.1715728752538057);
\draw [line width=1.6pt,color=ffqqqq] (13.,-0.1715728752538057)-- (11.5857864376269,-1.5857864376269015);
\draw [line width=1.6pt,color=ffqqqq] (11.5857864376269,-1.5857864376269015)-- (11.585786437626899,-3.585786437626904);
\draw [line width=1.6pt,color=qqffff] (11.585786437626899,-3.585786437626904)-- (13.,-5.);
\draw [line width=1.6pt,color=ffqqqq] (26.071067811865518,-1.5857864376268749)-- (24.656854249492298,-0.17157287525375875);
\draw [line width=1.6pt,color=ffqqqq] (24.656854249492298,-0.17157287525375875)-- (22.656854249492305,-0.17157287525388165);
\draw [line width=1.6pt,color=ffqqqq] (22.656854249492305,-0.17157287525388165)-- (21.242640687119305,-1.5857864376270585);
\draw [line width=0.0pt,color=ffqqqq] (21.242640687119305,-1.5857864376270585)-- (21.24264068711943,-3.585786437627049);
\draw [line width=1.6pt,color=qqffff] (21.24264068711943,-3.585786437627049)-- (22.656854249492426,-5.);
\draw [line width=1.6pt,color=qqqqff] (22.656854249492426,-5.)-- (24.656854249492444,-5.);
\draw [line width=1.6pt,color=ffxfqq] (16.414213562373057,6.071067811865503)-- (17.828427124746117,4.656854249492467);
\draw [line width=1.6pt,color=ttzzqq] (19.828427124746085,4.6568542494924685)-- (21.242640687119156,6.071067811865542);
\draw [line width=1.6pt,color=ffqqqq] (21.242640687119156,6.071067811865542)-- (21.242640687119156,8.07106781186551);
\draw [line width=1.6pt,color=ttzzqq] (21.242640687119156,8.07106781186551)-- (19.82842712474608,9.485281374238582);
\draw [line width=1.6pt,color=ttzzqq] (19.82842712474608,9.485281374238582)-- (17.828427124746113,9.48528137423858);
\draw [line width=1.6pt,color=ttzzqq] (17.828427124746113,9.48528137423858)-- (16.414213562373043,8.071067811865507);
\draw [line width=1.6pt,color=ttzzqq] (16.414213562373043,8.071067811865507)-- (16.414213562373057,6.071067811865503);
\draw [line width=1.6pt,color=ffqqqq] (21.24264068711943,-3.585786437627049)-- (21.242640687119305,-1.5857864376270585);
\draw [line width=1.6pt,color=ttzzqq] (21.242640687119305,-1.5857864376270585)-- (19.8284271247462,-0.17157287525380616);
\draw [line width=1.6pt,color=qqzzqq] (19.8284271247462,-0.17157287525380616)-- (17.828427124746195,-0.17157287525380438);
\draw [line width=1.6pt,color=ttzzqq] (17.828427124746195,-0.17157287525380438)-- (16.414213562373096,-1.5857864376269026);
\draw [line width=1.6pt,color=ttzzqq] (16.414213562373096,-1.5857864376269026)-- (16.414213562373096,-3.585786437626905);
\draw [line width=1.6pt,color=ffxfqq] (16.414213562373096,-3.585786437626905)-- (17.82842712474619,-5.);
\draw [line width=1.6pt,color=ffffqq] (17.82842712474619,-5.)-- (19.828427124746195,-5.);
\draw [line width=1.6pt,color=qqffff] (19.828427124746195,-5.)-- (21.24264068711943,-3.585786437627049);
\draw [line width=1.6pt,color=ttzzqq] (30.899494936611788,-1.5857864376267878)-- (29.485281374238532,-0.17157287525359122);
\draw [line width=1.6pt,color=ttzzqq] (29.485281374238532,-0.17157287525359122)-- (27.48528137423853,-0.17157287525356602);
\draw [line width=1.6pt,color=ttzzqq] (27.48528137423853,-0.17157287525356602)-- (26.071067811865518,-1.5857864376268749);
\draw [line width=1.6pt,color=ffffqq] (27.485281374238706,-5.)-- (29.485281374238745,-5.);
\draw [line width=1.6pt,color=qqffff] (29.485281374238745,-5.)-- (30.899494936611834,-3.585786437626825);
\draw [line width=1.6pt,color=ffqqqq] (30.899494936611834,-3.585786437626825)-- (30.899494936611788,-1.5857864376267878);
\draw [line width=1.6pt,color=ffffqq] (16.414213562373106,1.2426406871193165)-- (16.414213562373085,3.2426406871193167);
\draw [line width=0.0pt,color=qqqqff] (16.414213562373085,3.2426406871193167)-- (15.,4.6568542494923975);
\draw [line width=0.0pt,color=qqqqff] (15.,4.6568542494923975)-- (13.,4.656854249492376);
\draw [line width=0.0pt,color=qqqqff] (13.,4.656854249492376)-- (11.585786437626894,3.2426406871192666);
\draw [line width=0.0pt,color=qqqqff] (11.585786437626894,3.2426406871192666)-- (11.585786437626915,1.2426406871192657);
\draw [line width=0.0pt,color=qqqqff] (11.585786437626915,1.2426406871192657)-- (13.,-0.1715728752538057);
\draw [line width=0.0pt,color=qqqqff] (13.,-0.1715728752538057)-- (15.,-0.17157287525380602);
\draw [line width=0.0pt,color=qqqqff] (15.,-0.17157287525380602)-- (16.414213562373106,1.2426406871193165);
\draw [line width=0.0pt,color=qqqqff] (21.242640687119287,1.2426406871193783)-- (21.24264068711926,3.242640687119373);
\draw [line width=0.0pt,color=ffffqq] (21.24264068711926,3.242640687119373)-- (19.828427124746085,4.6568542494924685);
\draw [line width=0.0pt,color=ffffqq] (17.828427124746117,4.656854249492467)-- (16.414213562373085,3.2426406871193167);
\draw [line width=0.0pt,color=ffffqq] (16.414213562373085,3.2426406871193167)-- (16.414213562373106,1.2426406871193165);
\draw [line width=0.0pt,color=ffffqq] (16.414213562373106,1.2426406871193165)-- (17.828427124746195,-0.17157287525380438);
\draw [line width=1.6pt,color=qqzzqq] (17.828427124746195,-0.17157287525380438)-- (19.8284271247462,-0.17157287525380616);
\draw [line width=0.0pt,color=ffffqq] (19.8284271247462,-0.17157287525380616)-- (21.242640687119287,1.2426406871193783);
\draw [line width=1.6pt,color=ffqqqq] (22.656854249492305,-0.17157287525388165)-- (24.656854249492298,-0.17157287525375875);
\draw [line width=1.6pt,color=qqffff] (24.656854249492298,-0.17157287525375875)-- (26.071067811865483,1.242640687119466);
\draw [line width=1.6pt,color=ffffqq] (26.071067811865483,1.242640687119466)-- (26.071067811865436,3.2426406871194624);
\draw [line width=1.6pt,color=qqqqff] (26.071067811865436,3.2426406871194624)-- (24.65685424949222,4.6568542494924845);
\draw [line width=0.0pt,color=qqqqff] (24.65685424949222,4.6568542494924845)-- (22.656854249492255,4.656854249492461);
\draw [line width=0.0pt,color=qqqqff] (22.656854249492255,4.656854249492461)-- (21.24264068711926,3.242640687119373);
\draw [line width=1.6pt,color=qqqqff] (21.24264068711926,3.242640687119373)-- (21.242640687119287,1.2426406871193783);
\draw [line width=0.0pt,color=qqqqff] (21.242640687119287,1.2426406871193783)-- (22.656854249492305,-0.17157287525388165);
\draw [line width=0.0pt,color=ffxfqq] (26.071067811865273,6.071067811865571)-- (24.65685424949222,4.6568542494924845);
\draw [line width=1.6pt,color=qqqqff] (24.65685424949222,4.6568542494924845)-- (26.071067811865436,3.2426406871194624);
\draw [line width=1.6pt,color=ffffqq] (26.071067811865436,3.2426406871194624)-- (27.485281374238593,4.65685424949263);
\draw [line width=0.0pt,color=ffxfqq] (27.485281374238593,4.65685424949263)-- (26.071067811865273,6.071067811865571);
\draw [line width=0.0pt,color=ttffcc] (21.242640687119156,6.071067811865542)-- (19.828427124746085,4.6568542494924685);
\draw [line width=1.6pt,color=ffffqq] (19.828427124746085,4.6568542494924685)-- (21.24264068711926,3.242640687119373);
\draw [line width=1.6pt,color=qqqqff] (21.24264068711926,3.242640687119373)-- (22.656854249492255,4.656854249492461);
\draw [line width=0.0pt,color=ttffcc] (22.656854249492255,4.656854249492461)-- (21.242640687119156,6.071067811865542);
\draw [line width=1.6pt,color=qqwuqq] (17.82842712474618,-9.828427124746225)-- (19.82842712474619,-9.828427124746225);
\draw [line width=1.6pt,color=qqwuqq] (27.485281374238784,-9.828427124746222)-- (29.485281374238795,-9.828427124746213);
\draw [line width=1.6pt,color=qqqqff] (30.89949493661167,3.24264068711949)-- (30.899494936611646,1.2426406871194882);
\draw [line width=1.6pt,color=qqqqff] (30.899494936611887,-6.414213562373092)-- (30.899494936611895,-8.414213562373105);
\draw [line width=1.6pt,color=ffqqqq] (30.89949493661139,8.071067811865671)-- (30.899494936611468,6.0710678118656825);
\draw [line width=1.6pt,color=ffqqqq] (30.899494936611788,-1.5857864376267878)-- (30.899494936611834,-3.585786437626825);
\draw [line width=1.6pt,color=ttffcc] (16.414213562373106,1.2426406871193165)-- (17.828427124746195,-0.17157287525380438);
\draw [line width=1.6pt,color=ttffcc] (16.414213562373106,1.2426406871193165)-- (15.,-0.17157287525380602);
\draw [line width=1.6pt,color=qqffff] (19.828427124746085,4.6568542494924685)-- (21.242640687119156,6.071067811865542);
\draw [line width=1.6pt,color=qqffff] (21.242640687119156,6.071067811865542)-- (22.656854249492255,4.656854249492461);
\draw [line width=1.6pt,color=ffxfqq] (15.,4.6568542494923975)-- (16.414213562373057,6.071067811865503);
\draw [line width=1.6pt,color=ffxfqq] (16.414213562373057,6.071067811865503)-- (17.828427124746117,4.656854249492467);
\draw [line width=1.6pt,color=ffxfqq] (24.65685424949222,4.6568542494924845)-- (26.071067811865273,6.071067811865571);
\draw [line width=1.6pt,color=ffxfqq] (26.071067811865273,6.071067811865571)-- (27.485281374238593,4.65685424949263);
\draw [line width=1.6pt,color=ffxfqq] (19.8284271247462,-0.17157287525380616)-- (21.242640687119287,1.2426406871193783);
\draw [line width=1.6pt,color=ffxfqq] (21.242640687119287,1.2426406871193783)-- (22.656854249492305,-0.17157287525388165);
\draw [line width=1.6pt,color=qqffff] (29.485281374238593,4.6568542494926035)-- (30.899494936611468,6.0710678118656825);
\draw [line width=1.6pt,color=qqffff] (11.585786437626934,6.071067811865506)-- (13.,4.656854249492376);
\draw [line width=1.6pt,color=ffxfqq] (29.485281374238532,-0.17157287525359122)-- (30.899494936611646,1.2426406871194882);
\draw [line width=1.6pt,color=ffxfqq] (11.585786437626915,1.2426406871192657)-- (13.,-0.1715728752538057);
\draw [line width=1.6pt,color=qqffff] (26.071067811865483,1.242640687119466)-- (27.48528137423853,-0.17157287525356602);
\draw [line width=1.6pt,color=ffxfqq] (16.414213562373096,-3.585786437626905)-- (15.,-5.);
\draw [line width=1.6pt,color=ffxfqq] (26.071067811865536,-3.585786437626891)-- (24.656854249492444,-5.);
\draw [line width=1.6pt,color=ffxfqq] (26.071067811865536,-3.585786437626891)-- (27.485281374238706,-5.);
\draw [line width=1.6pt,color=ffxfqq] (29.485281374238795,-9.828427124746213)-- (30.899494936611895,-8.414213562373105);
\draw [line width=1.6pt,color=qqffff] (26.071067811865532,-8.414213562373087)-- (27.485281374238784,-9.828427124746222);
\draw [line width=1.6pt,color=qqffff] (24.656854249492447,-9.828427124746172)-- (26.071067811865532,-8.414213562373087);
\draw [line width=1.6pt,color=ffqqqq] (13.,-9.82842712474619)-- (15.,-9.828427124746192);
\draw [line width=1.6pt,color=qqffff] (15.,-9.828427124746192)-- (16.414213562373092,-8.414213562373098);
\draw [line width=1.6pt,color=qqffff] (16.414213562373092,-8.414213562373098)-- (17.82842712474618,-9.828427124746225);
\draw [line width=1.6pt,color=ffxfqq] (19.82842712474619,-9.828427124746225)-- (21.242640687119376,-8.414213562373087);
\draw [line width=1.6pt,color=ffxfqq] (21.242640687119376,-8.414213562373087)-- (22.65685424949246,-9.828427124746172);
\draw [line width=1.6pt,color=qqqqff] (11.585786437626904,-6.414213562373095)-- (13.,-5.);
\draw [line width=1.6pt,color=ffxfqq] (11.585786437626902,-8.414213562373094)-- (13.,-9.82842712474619);
\draw [line width=1.6pt,color=ffqqqq] (22.65685424949246,-9.828427124746172)-- (24.656854249492447,-9.828427124746172);
\draw [line width=1.6pt,color=zzttqq] (39.41421356237309,5.414213562373091)-- (38.,6.828427124746184);
\draw [line width=1.6pt,color=zzttqq] (38.,6.828427124746184)-- (36.,6.8284271247461845);
\draw [line width=1.6pt,color=zzttqq] (36.,6.8284271247461845)-- (34.58578643762691,5.41421356237309);
\draw [line width=1.6pt,color=zzttqq] (34.58578643762691,5.41421356237309)-- (34.58578643762691,3.414213562373092);
\draw [line width=1.0pt,dash pattern=on 2pt off 2pt,color=zzttqq] (34.58578643762691,3.414213562373092)-- (36.,2.);
\draw [line width=1.0pt,dash pattern=on 2pt off 2pt,color=zzttqq] (36.,2.)-- (38.,2.);
\draw [line width=1.0pt,dash pattern=on 2pt off 2pt,color=zzttqq] (38.,2.)-- (39.41421356237309,3.414213562373095);
\draw [line width=1.0pt,dash pattern=on 2pt off 2pt,color=zzttqq] (39.41421356237309,3.414213562373095)-- (39.41421356237309,5.414213562373091);
\draw [line width=1.6pt] (38.258333333333404,-1.96)-- (36.828333333333404,-0.5596853493024382);
\draw [line width=1.0pt,dash pattern=on 2pt off 2pt] (35.3983333333334,-2.02)-- (36.828333333333404,-3.420314650697562);
\draw [line width=1.0pt,dash pattern=on 2pt off 2pt] (36.828333333333404,-3.420314650697562)-- (38.258333333333404,-1.96);
\draw [line width=1.6pt] (36.828333333333404,-0.5596853493024382)-- (35.3983333333334,-2.02);
\draw (13.69833333333336,8.01) node[anchor=north west] {$1$};
\draw (23.418333333333376,7.95) node[anchor=north west] {$1$};
\draw (13.758333333333361,-1.65) node[anchor=north west] {$1$};
\draw (23.388333333333374,-1.62) node[anchor=north west] {$1$};
\draw (18.58833333333337,7.98) node[anchor=north west] {$2$};
\draw (28.15833333333338,7.95) node[anchor=north west] {$2$};
\draw (18.58833333333337,-1.68) node[anchor=north west] {$2$};
\draw (28.21833333333338,-1.71) node[anchor=north west] {$2$};
\draw (13.608333333333361,3.12) node[anchor=north west] {$3$};
\draw (23.388333333333374,3.21) node[anchor=north west] {$3$};
\draw (13.728333333333362,-6.54) node[anchor=north west] {$3$};
\draw (23.478333333333374,-6.48) node[anchor=north west] {$3$};
\draw (18.49833333333337,3.18) node[anchor=north west] {$4$};
\draw (28.12833333333338,3.18) node[anchor=north west] {$4$};
\draw (18.49833333333337,-6.39) node[anchor=north west] {$4$};
\draw (28.12833333333338,-6.48) node[anchor=north west] {$4$};
\draw (16.038333333333366,5.61) node[anchor=north west] {$5$};
\draw (25.698333333333377,5.61) node[anchor=north west] {$5$};
\draw (20.95833333333337,0.69) node[anchor=north west] {$5$};
\draw (16.008333333333365,-4.17) node[anchor=north west] {$5$};
\draw (25.75833333333338,-4.11) node[anchor=north west] {$5$};
\draw (20.86833333333337,5.58) node[anchor=north west] {$6$};
\draw (16.038333333333366,0.75) node[anchor=north west] {$6$};
\draw (25.728333333333378,0.75) node[anchor=north west] {$6$};
\draw (20.898333333333373,-4.11) node[anchor=north west] {$6$};
\draw (35.74833333333339,6.33) node[anchor=north west] {$\frac{1}{2}C$};
\begin{scriptsize}
\draw [fill=qqqqff] (13.,-5.) circle (5pt);
\draw [fill=qqqqff] (15.,-5.) circle (5pt);
\draw [fill=ffxfqq] (16.414213562373096,-3.585786437626905) circle (5pt);
\draw [fill=ttzzqq] (16.414213562373096,-1.5857864376269026) circle (5pt);
\draw [fill=qqffff] (16.414213562373106,1.2426406871193165) circle (5pt);
\draw [fill=ffffqq] (16.414213562373085,3.2426406871193167) circle (5pt);
\draw [fill=qqqqff] (15.,4.6568542494923975) circle (5pt);
\draw [fill=qqqqff] (13.,4.656854249492376) circle (5pt);
\draw [fill=qqqqff] (11.585786437626894,3.2426406871192666) circle (5pt);
\draw [fill=ffxfqq] (11.585786437626915,1.2426406871192657) circle (5pt);
\draw [fill=ffxfqq] (16.414213562373057,6.071067811865503) circle (5pt);
\draw [fill=ttzzqq] (16.414213562373043,8.071067811865507) circle (5pt);
\draw [fill=ffqqqq] (15.,9.485281374238548) circle (5pt);
\draw [fill=ffqqqq] (13.,9.485281374238548) circle (5pt);
\draw [fill=ffqqqq] (11.585786437626936,8.071067811865476) circle (5pt);
\draw [fill=qqffff] (11.585786437626934,6.071067811865506) circle (5pt);
\draw [fill=qqqqff] (13.,4.65685424949243) circle (5pt);
\draw [fill=qqqqff] (15.,4.6568542494924285) circle (5pt);
\draw [fill=ffqqqq] (15.,-0.17157287525380602) circle (5pt);
\draw [fill=ffqqqq] (13.,-0.1715728752538057) circle (5pt);
\draw [fill=ffqqqq] (11.5857864376269,-1.5857864376269015) circle (5pt);
\draw [fill=qqffff] (11.585786437626899,-3.585786437626904) circle (5pt);
\draw [fill=qqqqff] (13.,-5.) circle (5pt);
\draw [fill=qqqqff] (15.,-5.) circle (5pt);
\draw [fill=ffffqq] (17.82842712474619,-5.) circle (5pt);
\draw [fill=ffffqq] (19.828427124746195,-5.) circle (5pt);
\draw [fill=ttzzqq] (19.8284271247462,-0.17157287525380616) circle (5pt);
\draw [fill=ttzzqq] (17.828427124746195,-0.17157287525380438) circle (5pt);
\draw [fill=ttzzqq] (17.82842712474622,-0.17157287525376116) circle (5pt);
\draw [fill=ttzzqq] (19.828427124746216,-0.17157287525373377) circle (5pt);
\draw [fill=ffxfqq] (21.242640687119287,1.2426406871193783) circle (5pt);
\draw [fill=qqqqff] (21.24264068711926,3.242640687119373) circle (5pt);
\draw [fill=ffffqq] (17.828427124746117,4.656854249492467) circle (5pt);
\draw [fill=ffffqq] (19.828427124746085,4.6568542494924685) circle (5pt);
\draw [fill=qqffff] (21.242640687119156,6.071067811865542) circle (5pt);
\draw [fill=ffqqqq] (21.242640687119156,8.07106781186551) circle (5pt);
\draw [fill=ttzzqq] (19.82842712474608,9.485281374238582) circle (5pt);
\draw [fill=ttzzqq] (17.828427124746113,9.48528137423858) circle (5pt);
\draw [fill=qqqqff] (22.656854249492426,-5.) circle (5pt);
\draw [fill=qqqqff] (24.656854249492444,-5.) circle (5pt);
\draw [fill=ffxfqq] (26.071067811865536,-3.585786437626891) circle (5pt);
\draw [fill=ttzzqq] (26.071067811865518,-1.5857864376268749) circle (5pt);
\draw [fill=qqffff] (26.071067811865483,1.242640687119466) circle (5pt);
\draw [fill=ffffqq] (26.071067811865436,3.2426406871194624) circle (5pt);
\draw [fill=qqqqff] (22.656854249492255,4.656854249492461) circle (5pt);
\draw [fill=qqqqff] (24.65685424949222,4.6568542494924845) circle (5pt);
\draw [fill=ffxfqq] (26.071067811865273,6.071067811865571) circle (5pt);
\draw [fill=ttzzqq] (26.071067811865248,8.071067811865536) circle (5pt);
\draw [fill=ffqqqq] (24.656854249492163,9.485281374238589) circle (5pt);
\draw [fill=ffqqqq] (22.6568542494922,9.485281374238566) circle (5pt);
\draw [fill=ffqqqq] (24.656854249492298,-0.17157287525375875) circle (5pt);
\draw [fill=ffqqqq] (22.656854249492305,-0.17157287525388165) circle (5pt);
\draw [fill=ffqqqq] (21.242640687119305,-1.5857864376270585) circle (5pt);
\draw [fill=qqffff] (21.24264068711943,-3.585786437627049) circle (5pt);
\draw [fill=qqqqff] (22.656854249492604,-5.) circle (5pt);
\draw [fill=qqqqff] (24.656854249492593,-5.) circle (5pt);
\draw [fill=qqffff] (30.899494936611468,6.0710678118656825) circle (5pt);
\draw [fill=ffqqqq] (30.89949493661139,8.071067811865671) circle (5pt);
\draw [fill=ttzzqq] (29.485281374238248,9.485281374238703) circle (5pt);
\draw [fill=ttzzqq] (27.485281374238262,9.485281374238625) circle (5pt);
\draw [fill=ttzzqq] (27.48528137423853,-0.17157287525356602) circle (5pt);
\draw [fill=ttzzqq] (29.485281374238532,-0.17157287525359122) circle (5pt);
\draw [fill=ffxfqq] (30.899494936611646,1.2426406871194882) circle (5pt);
\draw [fill=qqqqff] (30.89949493661167,3.24264068711949) circle (5pt);
\draw [fill=ffffqq] (29.485281374238593,4.6568542494926035) circle (5pt);
\draw [fill=ffffqq] (27.485281374238593,4.65685424949263) circle (5pt);
\draw [fill=ffffqq] (27.485281374238706,-5.) circle (5pt);
\draw [fill=ffffqq] (29.485281374238745,-5.) circle (5pt);
\draw [fill=qqffff] (30.899494936611834,-3.585786437626825) circle (5pt);
\draw [fill=ffqqqq] (30.899494936611788,-1.5857864376267878) circle (5pt);
\draw [fill=ttzzqq] (29.485281374238635,-0.1715728752536987) circle (5pt);
\draw [fill=ttzzqq] (27.485281374238596,-0.17157287525374354) circle (5pt);
\draw [fill=qqqqff] (11.585786437626904,-6.414213562373095) circle (5pt);
\draw [fill=ffxfqq] (11.585786437626902,-8.414213562373094) circle (5pt);
\draw [fill=ffqqqq] (13.,-9.82842712474619) circle (5pt);
\draw [fill=ffqqqq] (15.,-9.828427124746192) circle (5pt);
\draw [fill=ffffqq] (16.414213562373092,-6.4142135623730985) circle (5pt);
\draw [fill=ffffqq] (16.414213562373078,-6.41421356237311) circle (5pt);
\draw [fill=qqffff] (16.414213562373078,-8.414213562373122) circle (5pt);
\draw [fill=qqzzqq] (17.82842712474618,-9.828427124746225) circle (5pt);
\draw [fill=qqzzqq] (19.82842712474619,-9.828427124746225) circle (5pt);
\draw [fill=qqqqff] (21.242640687119298,-6.414213562373112) circle (5pt);
\draw [fill=qqqqff] (21.242640687119376,-6.414213562373101) circle (5pt);
\draw [fill=ffxfqq] (21.242640687119376,-8.414213562373087) circle (5pt);
\draw [fill=ffqqqq] (22.65685424949246,-9.828427124746172) circle (5pt);
\draw [fill=ffqqqq] (24.656854249492447,-9.828427124746172) circle (5pt);
\draw [fill=ffffqq] (26.071067811865532,-6.414213562373101) circle (5pt);
\draw [fill=ffffqq] (26.071067811865667,-6.414213562373111) circle (5pt);
\draw [fill=ffffff] (26.071067811865674,-8.414213562373122) circle (5pt);
\draw [fill=qqzzqq] (27.485281374238784,-9.828427124746222) circle (5pt);
\draw [fill=qqzzqq] (29.485281374238795,-9.828427124746213) circle (5pt);
\draw [fill=ffxfqq] (30.899494936611895,-8.414213562373105) circle (5pt);
\draw [fill=qqqqff] (30.899494936611887,-6.414213562373092) circle (5pt);
\draw [fill=ffffff] (36.,2.) circle (5pt);
\draw [fill=ffffff] (38.,2.) circle (5pt);
\draw [fill=ffffff] (39.41421356237309,3.414213562373095) circle (5pt);
\draw [fill=ffffff] (39.41421356237309,5.414213562373091) circle (5pt);
\draw [fill=zzwwqq] (38.,6.828427124746184) circle (5pt);
\draw [fill=zzwwqq] (36.,6.8284271247461845) circle (5pt);
\draw [fill=zzwwqq] (34.58578643762691,5.41421356237309) circle (5pt);
\draw [fill=ffffff] (34.58578643762691,3.414213562373092) circle (5pt);
\draw [fill=ffffff] (35.3983333333334,-2.02) circle (5pt);
\draw [fill=ffffff] (38.258333333333404,-1.96) circle (5pt);
\draw [fill=uuuuuu] (36.828333333333404,-0.5596853493024382) circle (5pt);
\draw [fill=ffffff] (36.828333333333404,-3.420314650697562) circle (5pt);
\draw [fill=qqffff] (26.071067811865532,-8.414213562373087) circle (5pt);
\end{scriptsize}
\end{tikzpicture}
\caption{ A $6$-coloring of the Minkowski plane equipped with the regular octagon metric with no pair of points unit distance apart identically colored }
\label{fig3b}
\end{figure}

As noted by Chilakamarri \cite{chilakamarri}, when $C$ is a centrally symmetric octagon one can prove the better upper bound $\chi(\mathbf{R}^2, C)\le 6$ as there exists a $6$-coloring of the plane avoiding monochromatic unit distances. The tiling shown in Figure \ref{fig3b} illustrates the case when $C$ is a regular octagon; it uses four colors for the  octagonal tiles which are translates of $\frac{1}{2}C$ and two additional colors for the remaining square tiles. Note that for each tile, only half of the boundary points are colored the same as the interior.

\section{\bf $C=$ regular decagon }

Let $\theta=\pi/5$.

Consider the regular decagon with vertices located at $(2\cos{k\theta}, 2\sin{k\theta})$ with $k=0\ldots 9$. It is known that

\begin{equation*}
\cos{\theta}=\frac{1+\sqrt{5}}{4} \quad \text{and} \quad \sin{\theta}=\frac{\sqrt{10-2\sqrt{5}}}{4}.
\end{equation*}

\begin{figure}[htp]
\centering
\begin{tikzpicture}[line cap=round,line join=round,x=1.0cm,y=1.0cm,scale=0.46]
\clip(-9.86,-6.68) rectangle (11.1,6.98);
\draw (6.28,0.72) node[anchor=north west] {$\mathbf{a}_1(-2,2,0,0)$};
\draw (5.18,4.4) node[anchor=north west] {$\mathbf{a}_2(2,0,2,0)$};
\draw (1.9,7.08) node[anchor=north west] {$\mathbf{a}_3(3,-1,1,1)$};
\draw (-4.8,7.14) node[anchor=north west] {$\mathbf{a}_4(-3,1,1,1)$};
\draw (-9.82,4.5) node[anchor=north west] {$\mathbf{a}_5(-2,0,2,0)$};
\draw (-7.9,0.92) node[anchor=north west] {$-\mathbf{a}_1$};
\draw (-6.72,-3.06) node[anchor=north west] {$-\mathbf{a}_2$};
\draw (-3.08,-5.52) node[anchor=north west] {$-\mathbf{a}_3$};
\draw (1.06,-5.52) node[anchor=north west] {$-\mathbf{a}_4$};
\draw (4.86,-3.06) node[anchor=north west] {$-\mathbf{a}_5$};
\draw [line width=2.pt] (6.,0.)-- (4.854101966249685,3.526711513754839);
\draw [line width=2.pt] (4.854101966249685,3.526711513754839)-- (1.8541019662496847,5.706339097770922);
\draw [line width=2.pt] (1.8541019662496847,5.706339097770922)-- (-1.8541019662496845,5.706339097770923);
\draw [line width=2.pt] (-1.8541019662496845,5.706339097770923)-- (-4.854101966249685,3.5267115137548397);
\draw [line width=2.pt] (-4.854101966249685,3.5267115137548397)-- (-6.,0.);
\draw [line width=2.pt] (-6.,0.)-- (-4.854101966249685,-3.5267115137548384);
\draw [line width=2.pt] (-4.854101966249685,-3.5267115137548384)-- (-1.8541019662496847,-5.706339097770922);
\draw [line width=2.pt] (-1.8541019662496847,-5.706339097770922)-- (1.8541019662496845,-5.706339097770923);
\draw [line width=2.pt] (1.8541019662496845,-5.706339097770923)-- (4.854101966249685,-3.5267115137548397);
\draw [line width=2.pt] (4.854101966249685,-3.5267115137548397)-- (6.,0.);
\draw (-0.62,-0.04) node[anchor=north west] {$\mathbf{o}(0,0,0,0)$};
\begin{scriptsize}
\draw [fill=ffffff] (0.,0.) circle (2.5pt);
\draw [fill=ffqqqq] (6.,0.) circle (3.5pt);
\draw [fill=ffqqqq] (4.854101966249685,3.526711513754839) circle (3.5pt);
\draw [fill=ffqqqq] (1.8541019662496847,5.706339097770922) circle (3.5pt);
\draw [fill=ffqqqq] (-1.8541019662496845,5.706339097770923) circle (3.5pt);
\draw [fill=ffqqqq] (-4.854101966249685,3.5267115137548397) circle (3.5pt);
\draw [fill=ffqqqq] (-6.,0.) circle (3.5pt);
\draw [fill=ffqqqq] (-4.854101966249685,-3.5267115137548384) circle (3.5pt);
\draw [fill=ffqqqq] (-1.8541019662496847,-5.706339097770922) circle (3.5pt);
\draw [fill=ffqqqq] (1.8541019662496845,-5.706339097770923) circle (3.5pt);
\draw [fill=ffqqqq] (4.854101966249685,-3.5267115137548397) circle (3.5pt);
\end{scriptsize}
\end{tikzpicture}
\caption{ A regular decagon.}
\label{fig4}
\end{figure}
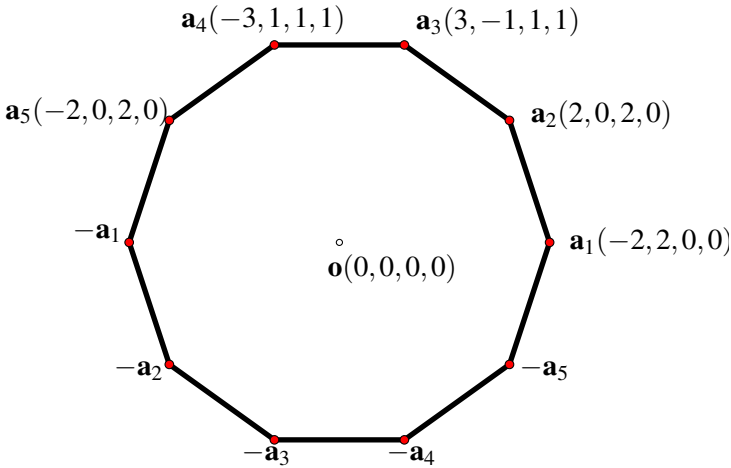

As we are trying to circumvent the nuisance of having to  manipulate (nested) radicals, we will only work with points whose coordinates are of the form
\begin{equation*}
\left((a+b\sqrt{5})\cos{\theta},  (c+d\sqrt{5})\sin{\theta}\right) \quad \text{where}\,\, a, b, c, d \quad \text{are integers}.
\end{equation*}
We will then write the coordinates of such a point as $(a,b,c,d)$. With this particular convention, the vertices of the decagon have particularly simple expressions. For example,
\begin{align*}
\mathbf{a}_1=&(2,0)=\left(\left(-2+2\sqrt{5}\right)\cos{\theta}, \left(0+0\sqrt{5}\right)\sin{\theta}\right) =(-2,2,0,0),\\
\mathbf{a}_2=&\left(2\cos{\theta}, 2\sin{\theta}\right)=\left(\left(2+0\sqrt{5}\right)\cos{\theta}, \left(2+0\sqrt{5}\right)\sin{\theta}\right)=(2,0,2,0),\\
\mathbf{a}_3=&\left(2\cos{2\theta}, 2\sin{2\theta}\right)=\left(\left(3-\sqrt{5}\right)\cos{\theta}, \left(1+\sqrt{5}\right)\sin{\theta}\right)=(3,-1,1,1).
\end{align*}
The coordinates of the remaining vertices can then be easily obtained by symmetry - see Figure \ref{fig4}.

For a given point $\mathbf{p}_1=(a,b,c,d)$, let $\mathbf{p_2}$, $\mathbf{p_3}$, $\mathbf{p}_4$, and $\mathbf{p}_5$  be the images of $\mathbf{p}_1$ after  counterclockwise rotations of $\pi/5, 2\pi/5$, $3\pi/5$, and $4\pi/5$, respectively. It is straightforward to verify that
\begin{align*}
&\mathbf{p_2}=\left(\frac{a+5b+5c-15d}{4}, \frac{a+b-3c+5d}{4}, \frac{a+5b+c+5d}{4}, \frac{a+b+c+d}{4}\right),\notag\\
&\mathbf{p_3}= \left(\frac{-a+5b-5c+5d}{4}, \frac{a-b+c-5d}{4}, \frac{3a+5b-c+5d}{4}, \frac{a+3b+c-d}{4}\right),\notag\\
&\mathbf{p_4}= \left(\frac{a-5b-5c+5d}{4}, \frac{a+b+c-5d}{4}, \frac{3a+5b+c-5d}{4}, \frac{a+3b-c+d}{4}\right),\\
&\mathbf{p_5}= \left(\frac{-a-5b+5c-15d}{4}, \frac{-a-b-3c+5d}{4}, \frac{a+5b-c-5d}{4}, \frac{a+b-c-d}{4}\right).\notag\\
\end{align*}

As in the octagon case, given a point $\mathbf{p}_1$ we denote by $\langle \mathbf{p}_1\rangle$ the set of vertices of the regular decagon centered at the origin one of whose vertices is $\mathbf{p}_1$.
\begin{equation}\label{decagonorbit}
\langle \mathbf{p}_1\rangle :=\{\mathbf{p}_1, \mathbf{p}_2, \mathbf{p}_3, \mathbf{p}_4, \mathbf{p}_5, -\mathbf{p}_1, -\mathbf{p}_2, -\mathbf{p}_3, -\mathbf{p}_4, -\mathbf{p}_5\}
\end{equation}
We refer to $\langle \mathbf{p}_1\rangle$ as the \emph{decagonal orbit} of $\mathbf{p}_1$.
While it is not difficult to find a unit-distance embedding of the Moser spindle in the regular decagon metric, we skip that step in favor of directly constructing a 5-chromatic unit distance graph. Rather surprisingly, we use fewer unit vectors than in the octagon case; moreover, the resulting graph does not contain 4-chromatic subgraphs with fewer that 10 vertices.

As in the octagon case, we present a set of \emph{generating} unit vectors whose 2-fold Minkowski sum will serve as a potential $5$-chromatic unit-distance graph, and a set of \emph{accidental} unit vectors which create additional edges within this graph. Since the total number of needed vectors is small, we introduce them all at once.

\begin{figure}[htp]
\centering
\begin{tikzpicture}[line cap=round,line join=round,x=1.0cm,y=1.0cm,scale=0.45]
\clip(-4.6525,-1.285) rectangle (23,10.835);
\draw (16.2675,0.815) node[anchor=north west] {$\mathbf{a}_1(-2,2,0,0)$};
\draw (13.1875,10.695) node[anchor=north west] {$\mathbf{a}_2(2,0,2,0)$};
\draw (0.8475,10.975) node[anchor=north west] {$\mathbf{a}_3$};
\draw (-1.7525,11.015) node[anchor=north west] {$\mathbf{a}_4$};
\draw (-3.6325,9.695) node[anchor=north west] {$\mathbf{a}_5$};
\draw (-4.6525,7.855) node[anchor=north west] {$-\mathbf{a}_1$};
\draw (-4.0725,5.835) node[anchor=north west] {$-\mathbf{a}_2$};
\draw (-2.1125,4.355) node[anchor=north west] {$-\mathbf{a}_3$};
\draw (0.5475,4.375) node[anchor=north west] {$-\mathbf{a}_4$};
\draw (2.4875,5.795) node[anchor=north west] {$-\mathbf{a}_5$};
\draw [line width=1.pt] (3.,7.)-- (2.4270509831248424,8.76335575687742);
\draw [line width=1.pt] (2.4270509831248424,8.76335575687742)-- (0.9270509831248419,9.853169548885461);
\draw [line width=1.pt] (0.9270509831248419,9.853169548885461)-- (-0.9270509831248426,9.853169548885461);
\draw [line width=1.pt] (-0.9270509831248426,9.853169548885461)-- (-2.4270509831248424,8.76335575687742);
\draw [line width=1.pt] (-2.4270509831248424,8.76335575687742)-- (-3.,7.);
\draw [line width=1.pt] (-3.,7.)-- (-2.427050983124843,5.23664424312258);
\draw [line width=1.pt] (-2.427050983124843,5.23664424312258)-- (-0.9270509831248424,4.146830451114538);
\draw [line width=1.pt] (-0.9270509831248424,4.146830451114538)-- (0.9270509831248429,4.146830451114538);
\draw [line width=1.pt] (0.9270509831248429,4.146830451114538)-- (2.4270509831248437,5.23664424312258);
\draw [line width=1.pt] (2.4270509831248437,5.23664424312258)-- (3.,7.);
\draw (-0.8325,7.115) node[anchor=north west] {$\mathbf{o}$};
\draw (3.1875,7.695) node[anchor=north west] {$\mathbf{a}_1$};
\draw (2.4275,9.795) node[anchor=north west] {$\mathbf{a}_2$};
\draw [line width=1pt,dash pattern=on 5pt off 5pt] (0.,0.)-- (12.94427190999916,9.40456403667957);
\draw [line width=1pt,dash pattern=on 5pt off 5pt] (0.,0.)-- (16.,0.);
\draw [line width=1pt] (16.,0.)-- (12.94427190999916,9.40456403667957);
\draw [line width=1pt,dash pattern=on 5pt off 5pt] (0.,7.)-- (3.,7.);
\draw [line width=1pt,dash pattern=on 5pt off 5pt] (0.,7.)-- (2.4270509831248424,8.76335575687742);
\draw (-2.5725,0.035) node[anchor=north west] {$\mathbf{o}(0,0,0,0)$};
\draw (14.9475,4.435) node[anchor=north west] {$\mathbf{b}_1(9, -3, 3, -1)$};
\draw (14.3475,6.715) node[anchor=north west] {$\mathbf{c}_1(-9,5,-1,1)$};
\draw [line width=1pt,dash pattern=on 5pt off 5pt] (0.,0.)-- (14.832815729997474,3.5922238126366826);
\draw [line width=1pt,dash pattern=on 5pt off 5pt] (0.,0.)-- (14.111456180001687,5.812340224042888);
\draw (15.5075,2.855) node[anchor=north west] {$\mathbf{d}_1(-20,10,-4,2)$};
\draw (13.7075,8.395) node[anchor=north west] {$\mathbf{e}_1(20,-8,6,-2)$};
\draw [line width=1pt,dotted] (0.,0.)-- (13.665631459994948,7.184447625273365);
\draw [line width=1pt,dotted] (0.,0.)-- (15.278640450004213,2.220116411406205);
\begin{scriptsize}
\draw [fill=ffffff] (0.,7.) circle (4.5pt);
\draw [fill=ffqqqq] (3.,7.) circle (4.5pt);
\draw [fill=ffqqqq] (2.4270509831248424,8.76335575687742) circle (4.5pt);
\draw [fill=ffqqqq] (0.9270509831248419,9.853169548885461) circle (4.5pt);
\draw [fill=ffqqqq] (-0.9270509831248426,9.853169548885461) circle (4.5pt);
\draw [fill=ffqqqq] (-2.4270509831248424,8.76335575687742) circle (4.5pt);
\draw [fill=ffqqqq] (-3.,7.) circle (4.5pt);
\draw [fill=ffqqqq] (-2.427050983124843,5.23664424312258) circle (4.5pt);
\draw [fill=ffqqqq] (-0.9270509831248424,4.146830451114538) circle (4.5pt);
\draw [fill=ffqqqq] (0.9270509831248429,4.146830451114538) circle (4.5pt);
\draw [fill=ffqqqq] (2.4270509831248437,5.23664424312258) circle (4.5pt);
\draw [fill=ffffff] (0.,0.) circle (4.5pt);
\draw [fill=ffqqqq] (16.,0.) circle (4.5pt);
\draw [fill=ffqqqq] (12.94427190999916,9.40456403667957) circle (4.5pt);
\draw [fill=ududff] (14.832815729997474,3.5922238126366826) circle (4.5pt);
\draw [fill=ududff] (0.36067977499790516,46.498721792343105) circle (4.5pt);
\draw [fill=ududff] (14.111456180001687,5.812340224042888) circle (4.5pt);
\draw [fill=ffffff] (15.278640450004213,2.220116411406205) circle (4.5pt);
\draw [fill=ffffff] (13.665631459994948,7.184447625273365) circle (4.5pt);
\end{scriptsize}
\end{tikzpicture}
\caption{ Generating vectors $\mathbf{a}_1, \mathbf{b}_1, \mathbf{c}_1$ in dashed lines; accidental vectors $\mathbf{d}_1, \mathbf{e}_1$ in dotted lines.}
\label{fig5}
\end{figure}
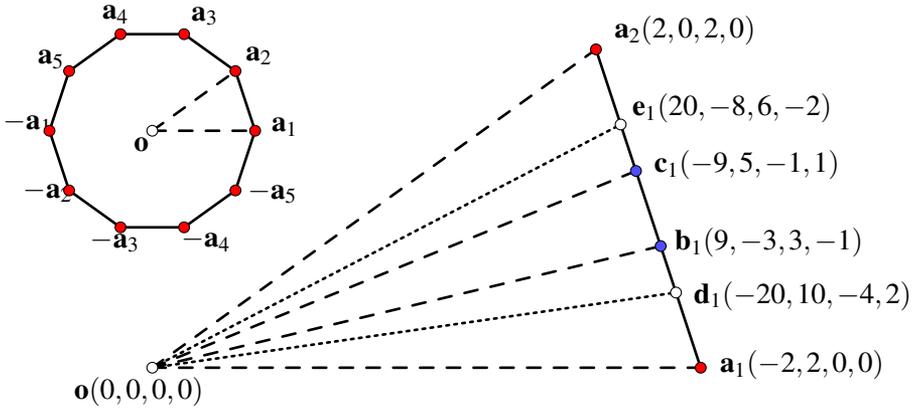

Let $\mathbf{b}_1, \mathbf{c}_1$ be the points which divide the side $\mathbf{a}_1\mathbf{a}_2$ in the ratios
$(3-\sqrt{5})/2$ and $(\sqrt{5}-1)/2$, respectively. Also, let $\mathbf{d}_1$ and $\mathbf{e}_1$ be the points which divide the side $\mathbf{a}_1\mathbf{a}_2$ in the ratios $\sqrt{5}-2$ and $3-\sqrt{5}$, respectively. See Figure \ref{fig5}.

It is easy to check that
\begin{align}\label{bcdedecagon}
&\mathbf{b}_1=(9,-3,3,-1), \mathbf{c}_1=(-9,5,-1,1),\notag\\
&\mathbf{d}_1=(-20,10,-4,2), \mathbf{e}_1=(20,-8,6,-2).
\end{align}
For verification purposes we write down explicitly the decagonal orbits of $\mathbf{a}_1, \mathbf{b}_1$, $\mathbf{c}_1$,
$\mathbf{d}_1$, and $\mathbf{e}_1$. In order to save space, only the first five vertices are listed in every orbit.
The remaining five vertices are obtained by taking opposites.
\small{
\begin{align*}
&\langle \mathbf{a}_1\rangle=\{(-2, 2, 0, 0), (2, 0, 2, 0), (3, -1, 1, 1), (-3, 1, 1, 1), (-2, 0, 2, 0),\ldots\},\\
&\langle \mathbf{b}_1\rangle=\{(9, -3, 3, -1), (6, -2, -2, 2), (-11, 5, 1, 1), (1, -1, 5, -1), (9, -5, -1, 1),\ldots\},\\
&\langle \mathbf{c}_1\rangle=\{(-9, 5, -1, 1), (-1, 1, 5, -1), (11, -5, 1, 1), (-6, 2, -2, 2), (-9, 3, 3, -1),\ldots\},\\
&\langle \mathbf{d}_1\rangle=\{(-20, 10, -4, 2), (-5, 3, 9, -3), (25, -11, 1, 1), (-10, 4, -6, 4), (-20, 8, 6, -2),\ldots\},\\
&\langle \mathbf{e}_1\rangle=\{(20, -8, 6, -2), (10, -4, -6, 4), (-25, 11, 1, 1), (5, -3, 9, -3), (20, -10, -4, 2),\ldots\}.
\end{align*}
}
Let $U_{30} = \langle \mathbf{a}_1\rangle \cup \langle \mathbf{b}_1\rangle \cup \langle \mathbf{c}_1\rangle$ be the promised set of generating unit vectors and let $V_{20}=\langle \mathbf{d}_1\rangle \cup \langle \mathbf{e}_1\rangle$ be the set of accidental unit vectors. The full set of unit vectors $W_{50}=U_{30}\cup V_{20}$ is all that is needed for constructing a $5$-chromatic graph.

Consider the $2$-fold Minkowski sum
\begin{equation}
U_{30}+U_{30} =\{ \mathbf{x}+\mathbf{y} \,\, | \,\, \mathbf{x}\in U_{30}, \mathbf{y} \in U_{30}\}.
\end{equation}

The resulting unit-distance graph with vertex set $U_{30}+U_{30}$ has $421$ vertices and $2640$ edges in $U_{30}$. Based solely on this edge set, the graph has chromatic number $4$. However, there are $500$ accidental pairs which determine vectors from $V_{20}$. The graph has now $3140$ edges and is $5$-chromatic.

\begin{figure}[ht]
\centering
\begin{tikzpicture}[line width=1pt,scale=1.5]
\tikzstyle{every node}=[draw=black,fill=yellow!50!white,thick,
  shape=circle,minimum height=0.2cm,inner sep=1];
\tikzstyle{every path}=[draw=white!10!black];

\input{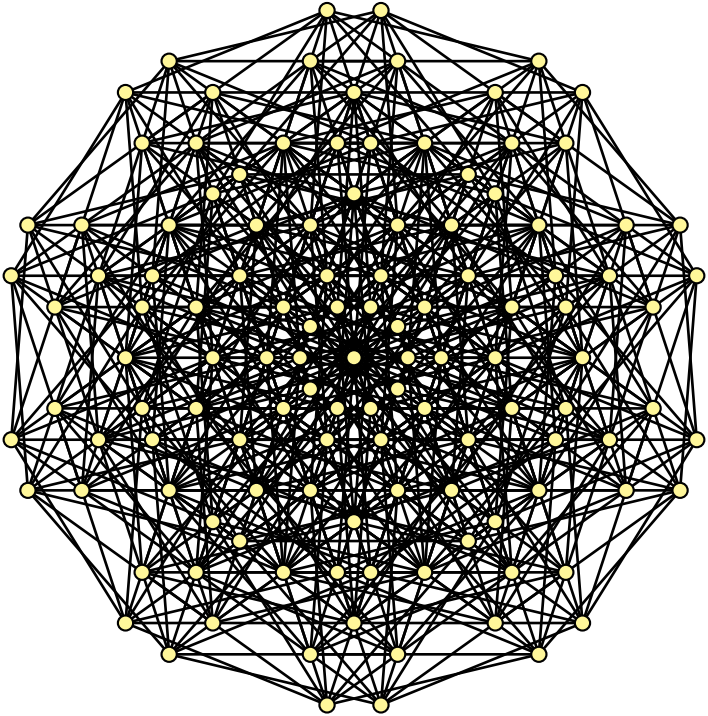}

\end{tikzpicture}
\caption{The expected edges in the $5$-chromatic graph $G_{121}$ for the regular decagon metric}
\label{decexpected}
\end{figure}

Again, one may inquire about a smaller $5$-chromatic graph. We found one such graph of order $121$.
\begin{thm}
Let $S$ consist of the following $13$ points in the Minkowski plane equipped with the regular decagon metric:
\begin{align*}
S:=\{& (0,0,0,0), (-9, 5, -1, 1), (-8, 4, 0, 0), (-6, 4, 2, 0), (-5, 3, 3, -1),\\
 & (-3, 3, -1, 1), (-2, 2, 0, 0), (1, 1, 3, -1), (5, -1, -1, 1) ,\\
 &(6, -2, 0, 0), (8, -2, 2, 0), (9, -3, 3, -1), (14, -6, 0, 0)\}
\end{align*}
Then the graph $G_{121}$, whose vertex set is $\bigcup_{i=1}^{13} \langle S_i \rangle$, has $121$ vertices, $680$ edges in $W_{50}$ and has chromatic number $5$. Hence, $\chi(\mathbf{R}^2, C)\ge 5$ when $C$ is a regular decagon.
\end{thm}
\begin{proof}
For verification purposes we provide the edge distribution over the unit vector orbits $\langle\mathbf{a}_1\rangle, \langle\mathbf{b}_1\rangle, \langle\mathbf{c}_1\rangle, \langle\mathbf{d}_1\rangle, \langle\mathbf{e}_1\rangle$:
\begin{equation*}
300, 180, 180, 10, 10.
\end{equation*}
\end{proof}
The graph is depicted in Figure~\ref{decexpected}, where edges of $G_{121}$ corresponding
to vectors in $U_{30}$ are shown, and Figure~\ref{decunexpected}, where edges of $G_{121}$ corresponding
to vectors in $V_{20}$ are shown.
As mentioned earlier, the graph $G_{121}$ contains no Moser spindles. In fact, the smallest $4$-chromatic subgraph has 10 vertices.
The graph $G_{121}$ contains (at least) $31$ pairwise non-isomorphic $4$-chromatic subgraphs of order $10$.
One of the more symmetric examples in shown in Figure~\ref{sub001}.

\begin{figure}
\centering
\begin{tikzpicture}[line width=1pt,scale=1.5]
\tikzstyle{every node}=[draw=black,fill=yellow!50!white,thick,
  shape=circle,minimum height=0.2cm,inner sep=1];
\tikzstyle{every path}=[draw=white!10!black];

\input{decdata2}

\end{tikzpicture}
\caption{The accidental edges in the $5$-chromatic graph  $G_{121}$ for the regular decagon metric}
\label{decunexpected}
\end{figure}

\begin{figure}
\centering
\begin{tikzpicture}[line width=1pt,scale=1.5]
\tikzstyle{every node}=[draw=black,fill=yellow!50!white,thick,
  shape=circle,minimum height=0.2cm,inner sep=1];
\tikzstyle{every path}=[draw=white!10!black];

\draw[black] (1.763932,-0.726543) -- (3.000000,0.726543) {};
\draw[black] (1.763932,-0.726543) -- (1.618034,1.175571) {};
\draw[black] (1.763932,-0.726543) -- (0.000000,0.000000) {};
\draw[black] (3.000000,0.726543) -- (1.381966,1.902113) {};
\draw[black] (3.000000,0.726543) -- (1.236068,0.000000) {};
\draw[black] (-0.236068,3.077684) -- (-1.236068,1.453085) {};
\draw[black] (-0.236068,3.077684) -- (-0.381966,1.175571) {};
\draw[black] (-0.236068,3.077684) -- (1.381966,1.902113) {};
\draw[black] (-1.236068,1.453085) -- (0.000000,0.000000) {};
\draw[black] (-1.236068,1.453085) -- (0.618034,1.902113) {};
\draw[black] (1.618034,1.175571) -- (0.000000,0.000000) {};
\draw[black] (1.618034,1.175571) -- (-0.381966,1.175571) {};
\draw[black] (0.000000,0.000000) -- (0.618034,1.902113) {};
\draw[black] (-0.381966,1.175571) -- (1.381966,1.902113) {};
\draw[black] (-0.381966,1.175571) -- (1.236068,0.000000) {};
\draw[black] (1.381966,1.902113) -- (1.236068,0.000000) {};
\draw[black] (1.236068,0.000000) -- (0.618034,1.902113) {};
\node at (1.763932,-0.726543) {};
\node at (3.000000,0.726543) {};
\node at (-0.236068,3.077684) {};
\node at (-1.236068,1.453085) {};
\node at (1.618034,1.175571) {};
\node at (0.000000,0.000000) {};
\node at (-0.381966,1.175571) {};
\node at (1.381966,1.902113) {};
\node at (1.236068,0.000000) {};
\node at (0.618034,1.902113) {};

\end{tikzpicture}
\caption{A $4$-chromatic subgraph of order $10$ in the graph $G_{121}$}
\label{sub001}
\end{figure}
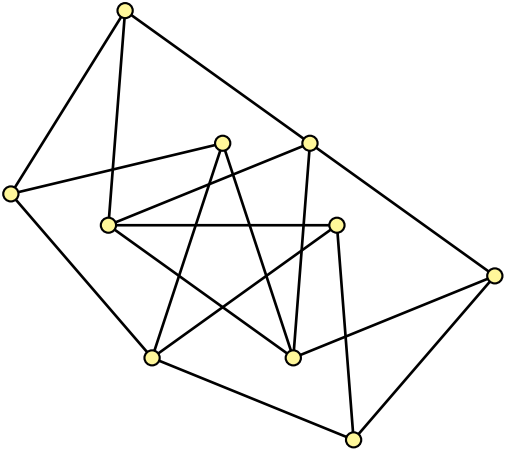

\newpage

\section{\bf $C=$ regular dodecagon}
Consider the regular dodecagon centered at the origin $\mathbf{o}(0,0)$ with vertices $\mathbf{a}_1=(12,0)$, $\mathbf{a}_2=(6\sqrt{3},6)$, $\mathbf{a}_3=(6,6\sqrt{3})$, $\mathbf{a}_4=(0,12)$, etc.

Throughout this section we will only work with points whose coordinates are of the form $(a+b\sqrt{3}, c+d\sqrt{3})$ with $a, b, c, d$ integers.
As in the previous cases, in order to avoid computations involving $\sqrt{3}$ we will encode such a point in the format $(a,b,c,d)$.
With this convention, the vertices of the regular octagon become
\begin{equation*}
\mathbf{a}_1=(12,0,0,0), \mathbf{a}_2=(0,6,6,0), \mathbf{a}_3=(6,0,0,6),\mathbf{a}_4=(0,0,12,0),\ldots
\end{equation*}

For a given point $\mathbf{p}_1=(a,b,c,d)$, let $\mathbf{p_2}$ and $\mathbf{p_3}$  be the images of $\mathbf{p}_1$ after  counterclockwise rotations of $\pi/3$ and $2\pi/3$, respectively.

It is easy to check that
\begin{equation}\label{hexagonorbit1}
\mathbf{p_2}=\left(\frac{a-3d}{2},\frac{b-c}{2},\frac{3b+c}{2}, \frac{a+d}{2}\right),
\mathbf{p_3}=\left(\frac{-a-3d}{2},\frac{-b-c}{2},\frac{3b-c}{2}, \frac{a-d}{2}\right)
\end{equation}

For a given point $\mathbf{p}_1$, we denote by $\langle \mathbf{p}_1\rangle$ the set of vertices of the regular \emph{hexagon} centered at the origin one of whose vertices is $\mathbf{p}_1$.
\begin{equation}\label{hexagonorbit2}
\langle \mathbf{p}_1\rangle :=\{\mathbf{p}_1, \mathbf{p}_2, \mathbf{p}_3,  -\mathbf{p}_1, -\mathbf{p}_2, -\mathbf{p}_3\}
\end{equation}
We refer to $\langle \mathbf{p}_1\rangle$ as the \emph{hexagonal orbit} of $\mathbf{p}_1$. Working with hexagonal orbit rather than with dodecagonal orbits allows us to better control the sets of generating and accidental vectors to be introduced below.

Let $\mathbf{b}_1$, $\mathbf{c}_1$, $\mathbf{d}_1$, and $\mathbf{e}_1$ be the points which divide the side $\mathbf{a}_1\mathbf{a}_2$ in the ratios $(2-\sqrt{3})/3$, $(5-2\sqrt{3})/6$, $(3-\sqrt{3})/3$, and $(4-\sqrt{3})/3$, respectively.

Quick computations reveal that
\begin{equation*}
\mathbf{b}_1=(-2, 8, 4, -2), \mathbf{c}_1=(-4, 9, 5, -2),\mathbf{d}_1=(-6, 10, 6, -2), \mathbf{e}_1=(-10, 12, 8, -2).
\end{equation*}
In addition, let $\mathbf{f}_1$, $\mathbf{g}_1$ be the points which divide the side $\mathbf{a}_2\mathbf{a}_3$ in the ratios $\sqrt{3}/6$ and $\sqrt{3}/3$, respectively. It is easy to check that
\begin{equation*}
\mathbf{f}_1=(-3, 7, 9, -1), \mathbf{g}_1=(-6, 8, 12, -2).
\end{equation*}
Define $U_{42}$, the set of generating vectors as
\begin{equation}\label{U42}
U_{42}:=\langle\mathbf{a}_2\rangle \cup \langle\mathbf{b}_1\rangle \cup \langle\mathbf{c}_1\rangle \cup \langle\mathbf{d}_1\rangle \cup \langle\mathbf{e}_1\rangle \cup \langle\mathbf{f}_1\rangle \cup \langle\mathbf{g}_1\rangle.
\end{equation}

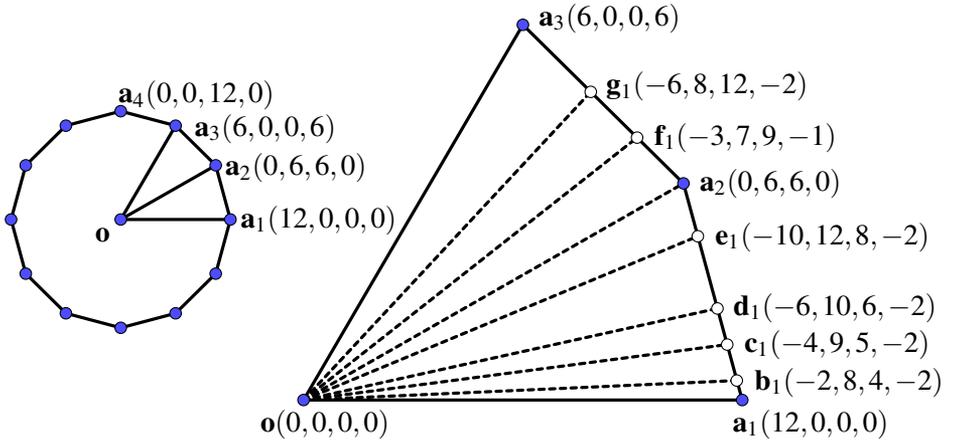
\begin{figure}[htp]
\centering
\begin{tikzpicture}[line cap=round,line join=round,x=1.0cm,y=1.0cm,scale=0.48]
\clip(-3.2,-6.12) rectangle (22.5,6.14);
\draw [line width=1.2pt] (3.,0.)-- (2.598076211353316,1.5);
\draw [line width=1.2pt] (1.5,2.5980762113533156)-- (2.598076211353316,1.5);
\draw [line width=1.22pt] (1.5,2.5980762113533156)-- (0.,3.);
\draw [line width=1.2pt] (0.,3.)-- (-1.5,2.5980762113533165);
\draw [line width=1.2pt] (-1.5,2.5980762113533165)-- (-2.598076211353315,1.5);
\draw [line width=1.2pt] (-2.598076211353315,1.5)-- (-3.,0.);
\draw [line width=1.2pt] (-3.,0.)-- (-2.598076211353317,-1.5);
\draw [line width=1.2pt] (-2.598076211353317,-1.5)-- (-1.5,-2.598076211353315);
\draw [line width=1.2pt] (-1.5,-2.598076211353315)-- (0.,-3.);
\draw [line width=1.2pt] (0.,-3.)-- (1.5,-2.5980762113533173);
\draw [line width=1.2pt] (1.5,-2.5980762113533173)-- (2.5980762113533147,-1.5);
\draw [line width=1.2pt] (2.5980762113533147,-1.5)-- (3.,0.);
\draw (-0.98,0.04) node[anchor=north west] {$\mathbf{o}$};
\draw (3.0,0.66) node[anchor=north west] {$\mathbf{a}_1(12,0,0,0)$};
\draw (2.58,2.12) node[anchor=north west] {$\mathbf{a}_2(0,6,6,0)$};
\draw (1.74,3.2) node[anchor=north west] {$\mathbf{a}_3(6,0,0,6)$};
\draw (-0.34,4.12) node[anchor=north west] {$\mathbf{a}_4(0,0,12,0)$};
\draw [line width=1.2pt] (5.,-5.)-- (17.,-5.);
\draw [line width=1.2pt,dash pattern=on 2pt off 2pt] (5.,-5.)-- (15.392304845413264,1.);
\draw [line width=1.2pt] (5.,-5.)-- (11.,5.392304845413262);
\draw [line width=1.2pt] (11.,5.392304845413262)-- (15.392304845413264,1.);
\draw [line width=1.2pt] (15.392304845413264,1.)-- (17.,-5.);
\draw (16.46,-4.98) node[anchor=north west] {$\mathbf{a}_1(12,0,0,0)$};
\draw (15.56,1.68) node[anchor=north west] {$\mathbf{a}_2(0,6,6,0)$};
\draw (11.16,6.26) node[anchor=north west] {$\mathbf{a}_3(6,0,0,6)$};
\draw (17.06,-3.82) node[anchor=north west] {$\mathbf{b}_1(-2,8,4,-2)$};
\draw (16.78,-2.78) node[anchor=north west] {$\mathbf{c}_1(-4,9,5,-2)$};
\draw (16.48,-1.76) node[anchor=north west] {$\mathbf{d}_1(-6,10,6,-2)$};
\draw (15.98,0.2) node[anchor=north west] {$\mathbf{e}_1(-10,12,8,-2)$};
\draw (14.32,2.96) node[anchor=north west] {$\mathbf{f}_1(-3,7,9,-1)$};
\draw (13.,4.38) node[anchor=north west] {$\mathbf{g}_1(-6,8,12,-2)$};
\draw (3.54,-5.) node[anchor=north west] {$\mathbf{o}(0,0,0,0)$};
\draw [line width=1.2pt] (0.,0.)-- (3.,0.);
\draw [line width=1.2pt] (0.,0.)-- (2.598076211353316,1.5);
\draw [line width=1.2pt] (0.,0.)-- (1.5,2.5980762113533156);
\draw [line width=1.2pt,dash pattern=on 2pt off 2pt] (5.,-5.)-- (16.85,-4.46);
\draw [line width=1.2pt,dash pattern=on 2pt off 2pt] (5.,-5.)-- (16.588,-3.464);
\draw [line width=1.2pt,dash pattern=on 2pt off 2pt] (5.,-5.)-- (16.3205,-2.4641);
\draw [line width=1.2pt,dash pattern=on 2pt off 2pt] (5.,-5.)-- (15.7846,-0.464);
\draw [line width=1.2pt,dash pattern=on 2pt off 2pt] (5.,-5.)-- (14.124,2.268);
\draw [line width=1.2pt,dash pattern=on 2pt off 2pt] (5.,-5.)-- (12.856,3.535);
\begin{scriptsize}
\draw [fill=ududff] (0.,0.) circle (4.5pt);
\draw [fill=ududff] (3.,0.) circle (4.5pt);
\draw [fill=ududff] (2.598076211353316,1.5) circle (4.5pt);
\draw [fill=ududff] (1.5,2.5980762113533156) circle (4.5pt);
\draw [fill=ududff] (0.,3.) circle (4.5pt);
\draw [fill=ududff] (-1.5,2.5980762113533165) circle (4.5pt);
\draw [fill=ududff] (-2.598076211353315,1.5) circle (4.5pt);
\draw [fill=ududff] (-3.,0.) circle (4.5pt);
\draw [fill=ududff] (-2.598076211353317,-1.5) circle (4.5pt);
\draw [fill=ududff] (-1.5,-2.598076211353315) circle (4.5pt);
\draw [fill=ududff] (0.,-3.) circle (4.5pt);
\draw [fill=ududff] (1.5,-2.5980762113533173) circle (4.5pt);
\draw [fill=ududff] (2.5980762113533147,-1.5) circle (4.5pt);
\draw [fill=ududff] (5.,-5.) circle (4.5pt);
\draw [fill=ududff] (17.,-5.) circle (4.5pt);
\draw [fill=ududff] (15.392304845413264,1.) circle (4.5pt);
\draw [fill=ududff] (11.,5.392304845413262) circle (4.5pt);
\draw [fill=ffffff] (16.85,-4.46) circle (4.5pt);
\draw [fill=ffffff] (16.588,-3.464) circle (4.5pt);
\draw [fill=ffffff] (16.3205,-2.4641) circle (4.5pt);
\draw [fill=ffffff] (15.7846,-0.464) circle (4.5pt);
\draw [fill=ffffff] (14.124,2.268) circle (4.5pt);
\draw [fill=ffffff] (12.856,3.535) circle (4.5pt);
\end{scriptsize}
\end{tikzpicture}
\caption{ Generating vectors $\mathbf{a}_2, \mathbf{b}_1, \mathbf{c}_1, \mathbf{d}_1, \mathbf{e}_1, \mathbf{f}_1, \mathbf{g}_1$ in dashed lines.}
\label{fig6}
\end{figure}

The explicit expressions of all vectors in $U_{42}$ can be easily obtained by applying \eqref{hexagonorbit1} to each of the relevant points - see Figure \ref{fig6}.

The graph with vertex set $U_{42}+U_{42}=\{\mathbf{x}+\mathbf{y}\,\,| \,\, \mathbf{x},\mathbf{y} \in U_{42}\}$ has 847 vertices, 4809 edges in $U_{42}$ and has chromatic number $4$. As in the previous cases, there exist many accidental unit vectors.

Let $\mathbf{h}_1$, $\mathbf{i}_1$, $\mathbf{j}_1$, $\mathbf{k}_1$, $\mathbf{l}_1$, $\mathbf{m}_1$, $\mathbf{n}_1$, $\mathbf{p}_1$, $\mathbf{q}_1$ be the points which divide the side $\mathbf{a}_1\mathbf{a}_2$ in the ratios $(7-4\sqrt{3})/6$, $(4-2\sqrt{3})/3$, $1/3$, $1/2$, $(5-2\sqrt{3})/3$, $(7-2\sqrt{3})/6$, $2/3$, $5/6$, and $(9-2\sqrt{3})/6$, respectively. Also, let $\mathbf{r}_1$ be the point which divides the side $\mathbf{a}_2\mathbf{a}_3$ in the ratio $\sqrt{3}/2$.

Straightforward computations show that
\begin{align*}
&\mathbf{h}_1=(-14,15,7,-4), \mathbf{i}_1=(-16,16,8,-4), \mathbf{j}_1=(8,2,2,0), \mathbf{k}_1=(6,3,3,0),\\
&\mathbf{l}_1=(-20,18,10,-4), \mathbf{m}_1=(-8,11,7,-2), \mathbf{n}_1=(4,4,4,0), \mathbf{p}_1=(2,5,5,0),\\
&\mathbf{q}_1=(-12,13,9,-2), \mathbf{r}_1=(-9,9,15,-3).
\end{align*}
Let $V_{60}$ be the set of accidental vectors defined as
\begin{equation*}
V_{60}=\langle\mathbf{h}_1\rangle \cup \langle\mathbf{i}_1\rangle \cup \langle\mathbf{j}_1\rangle \cup \langle\mathbf{k}_1\rangle \cup \langle\mathbf{l}_1\rangle \cup \langle\mathbf{m}_1\rangle \cup \langle\mathbf{n}_1\rangle \cup \langle\mathbf{p}_1\rangle \cup \langle\mathbf{q}_1\rangle \cup \langle\mathbf{r}_1\rangle.
\end{equation*}
It can be verified that the graph with vertex set $U_{42}+U_{42}$ has 1686 edges in $V_{60}$ besides the 4809 edges in $U_{42}$. The resulting unit distance graph has 847 vertices, 6495 edges, and has chromatic number 5. One may attempt to find a smaller $5$-chromatic graph: we found one such graph of order $295$.
\begin{thm}
Let $S$ consist of the following $50$ points in the Minkowski plane equipped with the regular dodecagon metric:
\begin{align*}
&(0, 0, 0, 0),(-12, 16, 0, 0), (-12, 16, 24, -4), (-12, 18, 6, 0), (-12, 18, 18, -4),\\
&(-12, 20, 0, 0), (-9, 14, 24, -5), (-9, 16, 6, -1), (-8, 10, -2, 4), (-8, 10, 10, -4),\\
& (-8, 12, 4, 0), (-6, 6, 6, -2), (-6, 7, 3, 0), (-6, 8, 0, 2), (-6, 8, 12, -6),\\
& (-6, 8, 12, -2), (-6, 9, 9, 0), (-6, 10, 6, -2), (-6, 10, 18, -2), (-6, 12, 0, 2),\\
& (-6, 12, 12, -2), (-6, 14, 6, -2), (-6, 14, 18, -2), (-6, 16, 12, -2), (-5, 7, 1, 1),\\
& (-5, 9, 7, -3), (-5, 15, 13, -3), (-4, 4, -4, 4), (-4, 6, 2, 0), (-4, 9, 17, -2),
\end{align*}

\begin{align*}
& (-4, 12, 8, 0), (-3, 6, 12, -3), (-3, 7, 9, -1), (-3, 8, 6, -3), (-3, 8, 6, 1),\\
& (-2, 4, 4, -2), (-2, 4, 4, 2), (-2, 8, 4, -2), (-2, 8, 16, -2), (0, 1, -3, 2),\\
& (0, 2, -6, 4), (0, 4, 0, 0), (0, 4, 12, -4), (0, 5, 9, -2), (0, 6, 6, 0), (0, 8, 24, -4),\\
& (2, 8, 8, -2), (2, 8, 20, -2), (4, 0, 4, 0), (6, -2, 6, -2).
\end{align*}

Then the graph $G_{295}$, whose vertex set is $\bigcup_{i=1}^{50} \langle S_i \rangle$ has $295$ vertices, $1644$ edges in $U_{42}\cup V_{60}$ and has chromatic number $5$. Hence, $\chi(\mathbf{R}^2, C)\ge 5$ when $C$ is a regular dodecagon.
\end{thm}
\begin{proof}
The edge distribution over the orbits of the generating vectors $\langle\mathbf{a}_2\rangle$, $\langle\mathbf{b}_1\rangle$, $\langle\mathbf{c}_1\rangle$, $\langle\mathbf{d}_1\rangle$, $\langle\mathbf{e}_1\rangle$, $\langle\mathbf{f}_1\rangle$, $\langle\mathbf{g}_1\rangle$ is 270, 168, 72, 222, 78, 162, 306. Similarly, the edge distribution over the orbits of the accidental vectors
$\langle\mathbf{h}_1\rangle$, $\langle\mathbf{i}_1\rangle$, $\langle\mathbf{j}_1\rangle$, $\langle\mathbf{k}_1\rangle$, $\langle\mathbf{l}_1\rangle$, $\langle\mathbf{m}_1\rangle$, $\langle\mathbf{n}_1\rangle$, $\langle\mathbf{p}_1\rangle$, $\langle\mathbf{q}_1\rangle$, $\langle\mathbf{r}_1\rangle$ is 6, 12, 42, 18, 0, 48, 60, 54, 48, 78.
\end{proof}

\begin{figure}[ht]
\centering
\begin{tikzpicture}[line width=1pt,scale=0.21]
\tikzstyle{every node}=[draw=black,fill=yellow!50!white,thick,
  shape=circle,minimum height=0.2cm,inner sep=1];
\tikzstyle{every path}=[draw=white!10!black];

\input{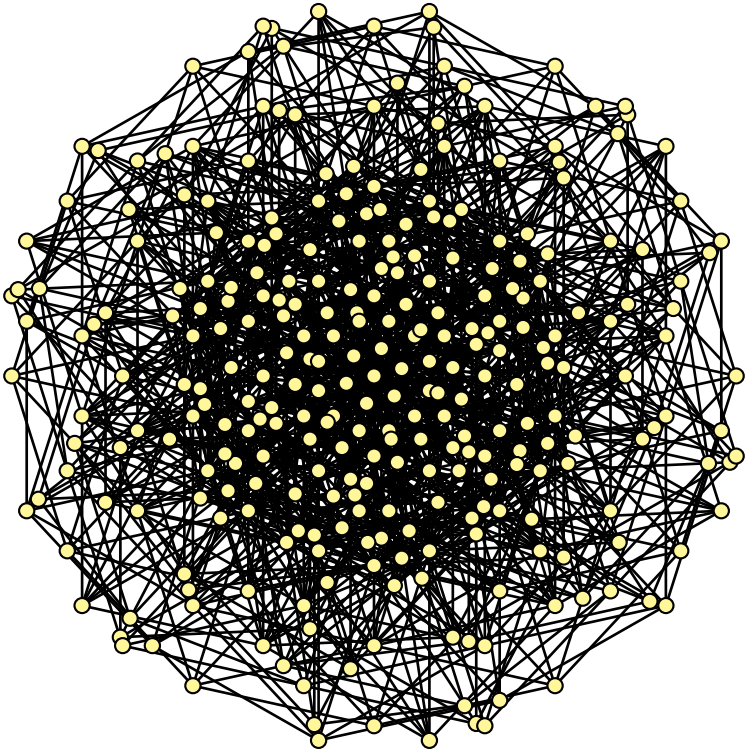}

\end{tikzpicture}
\caption{The expected edges in the graph $G_{295}$ for the dodecagon metric}
\label{dodexpected}
\end{figure}

\begin{figure}[ht]
\centering
\begin{tikzpicture}[line width=1pt,scale=0.21]
\tikzstyle{every node}=[draw=black,fill=yellow!50!white,thick,
  shape=circle,minimum height=0.2cm,inner sep=1];
\tikzstyle{every path}=[draw=white!10!black];

\input{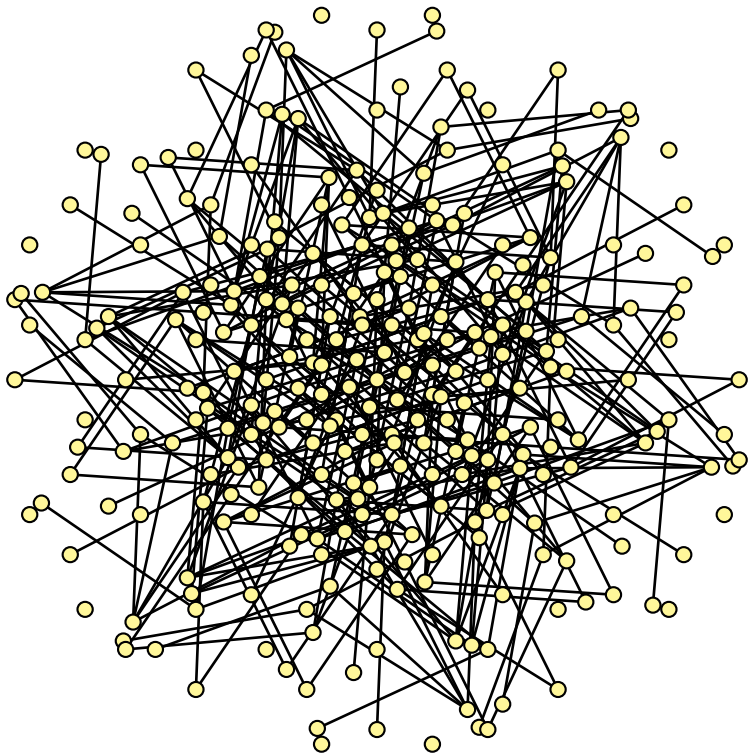}

\end{tikzpicture}
\caption{The accidental edges in the graph $G_{295}$ for the dodecagon metric}
\label{dodunexpected}
\end{figure}

We depict the graph $G_{295}$ in Figure~\ref{dodexpected} and
Figure~\ref{dodunexpected}, where the edges corresponding to vectors in
$U_{42}$ and $V_{60}$, respectively, are shown.

\section{\bf Open Problems}

In this paper we proved that the chromatic number of the Minkowski plane with the regular octagon, the regular decagon, or the regular dodecagon metric is $\ge 5$. One may wonder whether our approach can be extended to all regular $2n$-gons.

\begin{problem}
Is it true that $\chi(\mathbf{R}^2, C)\ge 5$ for every regular $2n$-gon with $n\ge 4$? More generally, is this true for every centrally symmetric
(not necessarily regular) $2n$-gon?
\end{problem}

\begin{problem}
Find a strictly convex, centrally symmetric curve $C$, different from an ellipse, for which  $\chi(\mathbf{R}^2, C)\ge 5$.
In particular, can one find some value $p>1, p\neq 2$, with $C:=\{(x,y)\,:\, |x|^p+|y|^p=1\}$ such that $\chi(\mathbf{R}^2, C)\ge 5$?
\end{problem}

\begin{problem}
Find a centrally symmetric convex curve $C$, other than the parallelogram or hexagon, for which  $\chi(\mathbf{R}^2, C)$ can be evaluated exactly.
\end{problem}

The most promising candidate seems to be the regular octagon as in this case the chromatic number of the Minkowski plane is either $5$ or $6$.
Small $4$-chromatic unit-distance graphs (in particular, the Moser spindle) were important in building $5$-chromatic unit-distance graphs in the Euclidean setting. It is reasonable to expect that small $5$-chromatic unit-distance graphs in some Minkowski plane would serve a similar purpose towards constructing a $6$-chromatic unit-distance graph.

\begin{problem}
Find a centrally symmetric convex curve $C$, for which there exists a $5$-chromatic graph in $(\mathbb{R}, C)$ of order $<100$.
\end{problem}


\begin{thebibliography}{99}


\bibitem{chilakamarri} K. B. Chilakamarri, Unit-distance graphs in Minkowski metric spaces, \emph{Geom. Dedicata} {\bf 37} (1991), no. 3, 345--356.

\bibitem{exooismailescuchi5}  G. Exoo and D. Ismailescu, The chromatic number of the plane is at least 5: a new proof, \emph{Discrete Comput. Geom.} {\bf 64} (2020), no. 1, 216--226.

\bibitem{degrey} A. D. N. J. de Grey, The chromatic number of the plane is at least $5$, \emph{Geombinatorics} {\bf 28} (2018), no. 1, 18--31.

\bibitem{hadwiger} H. Hadwiger,  \"{U}berdeckung des Euklidischen Raumes durch kongruente Mengen. (German), \emph{Portugaliae Math.} {\bf 4}, (1945). 238--242.

\bibitem{heule} M. J. H. Heule, Computing small unit-distance graphs with chromatic number 5, \emph{Geombinatorics} {\bf 28} (2018), no. 1, 32--50.

\bibitem{martini} H. Martini, K. J. Swanepoel, G. Weiss, The geometry of Minkowski spaces — a survey. I, \emph{Expo. Math.} {\bf 19} (2001), no. 2, 97--142.

\bibitem{MM} L. Moser, W. Moser, Solution to Problem $10$, \emph{Canadian Mathematical Bulletin} {\bf 4} (1961), 187--189.

\bibitem{parts1} J. Parts, Graph minimization, focusing on the example of 5-chromatic unit-distance graphs in the plane, \emph{Geombinatorics} {\bf 29} (2020), no. 3, 137--166.

\bibitem{parts2} J. Parts, The chromatic number of the plane is at least 5, a human-verifiable proof, \emph{Geombinatorics} {\bf 30} (2020), no. 2, 77--102.

\bibitem{polymath16} Polymath16, \url{https://asone.ai/polymath/index.php?title=Hadwiger-Nelson_problem.}



\bibitem{soifer} A. Soifer, {\it  The mathematical coloring book. Mathematics of coloring and the colorful life of its creators.} With forewords by Branko Gr\"{u}nbaum, Peter D. Johnson, Jr. and Cecil Rousseau. Springer, New York, 2009.

\end{thebibliography}
\end{document}